\documentclass{imsart}
\pubyear{0000}
\volume{00}
\issue{0}
\doi{0000}
\firstpage{1}
\lastpage{1}

\usepackage{amsthm}
\usepackage{amsmath}
\usepackage{natbib}
\usepackage[colorlinks,citecolor=blue,urlcolor=blue,filecolor=blue,backref=page]{hyperref}
\usepackage{graphicx}
\usepackage{multicol}
\usepackage{caption}
\usepackage{subcaption}

\usepackage{lineno,hyperref}
\usepackage{amssymb}

\usepackage[utf8]{inputenc}
\usepackage[english]{babel}

\startlocaldefs
\numberwithin{equation}{section}
\theoremstyle{plain}

\endlocaldefs

\theoremstyle{definition}

\begin{document}

\begin{frontmatter}
\title{Spatio-temporal stick-breaking process}
\runtitle{Spatio-temporal SB}

\begin{aug}
\author{\fnms{Clara} \snm{Grazian}\thanksref{addr1,addr2}\ead[label=e1]{clara.grazian@sydney.edu.au}},

\runauthor{Grazian, C.}

\address[addr1]{Carslaw Building, University of Sydney, Camperdown Campus, Sydney, Australia
    \printead{e1} 
}
\address[addr2]{ARC Training Centre in Data Analytics for Resources and Environments (DARE), Sydney, Australia
}
%

\end{aug}

\begin{abstract}
Dirichlet processes and their extensions have reached a great popularity in Bayesian nonparametric statistics. They have also been introduced for spatial and spatio-temporal data, as a tool to analyze and predict surfaces. 
A popular approach to Dirichlet processes in a spatial setting relies on a stick-breaking representation of the process, where the dependence over space is described in the definition of the stick-breaking probabilities. Extensions to include temporal dependence are still limited, however it is important, in particular for those phenomena which may change rapidly over time and space, with many local changes. In this work, we propose a Dirichlet process where the stick-breaking probabilities are defined to incorporate both spatial and temporal dependence. We will show that this approach is not a simple extension of available methodologies and can outperform available approaches in terms of prediction accuracy. An advantage of the method is that it offers a natural way to test for separability of the two components in the definition of the stick-breaking probabilities. 
\end{abstract}



\end{frontmatter}

\section{Introduction}

Spatial data is data that describes a surface, usually through variables on a two-dimensional plane, i.e. data that references a specific geographical area or location. The analysis and prediction of surfaces are essential to a diverse range of fields, including ecology, meteorology, sociology, image recognition, real estate, etc.

Spatial data is defined as a realisation of a stochastic process over space
$$
y(s) = \{y(s), s \in \mathcal{D}\}
$$
where $\mathcal{D} \subset \mathbb{R}^d$. The collected data comes as $y = \{y(s_1), \ldots, y(s_n)\}$, where $(s_1, \ldots, s_n)$ is a set of spatial locations, such that $s_i = (s_{i,1},s_{i,2})$ for $i=1,\ldots,n$; for example, $s_i$ may represent latitude and longitude. 
Depending on $\mathcal{D}$ being a continuous surface or a countable collection of $d$-dimensional spatial units, the problem can be specified as a spatially continuous or discrete random process, respectively \citep{gelfand2010handbook}. For example, data can be obtained by monitoring stations located in the set $(s_1, \ldots, s_n)$ of $n$ points. Alternatively, data may be observed in a set $(s_1, \ldots, s_n)$ of $n$ areas, defined, for example, by counties; we refer to this case as areal data. 

The first step in defining a spatial model within a Bayesian framework is to identify a probability distribution for the observed data. We usually select a distribution from the exponential family, indexed by a set of parameters $\theta$ accounting for spatial correlation; for instance
\begin{equation}
y(s) = \theta(s) + x(s)^T\beta + \varepsilon(s) 
\label{eq:spatmod}
\end{equation}
where $\theta(s) \in \Theta$ is the spatial effect, $x(s)^T$ is a set of $p$ covariates which may be dependent on the location $s$ and $\beta$ is a column-vector of $p$ regression coefficients; finally, $\varepsilon(s)$ is a normal error such that $\varepsilon(s) \sim \mathcal{N}(0,\sigma^2_{\varepsilon})$ and errors are considered independent and identically distributed (i.i.d.) over space.  

A common approach to the problem of modelling surfaces is kriging \citep{stein2012interpolation}, which interpolates and smooths spatial data; in kriging, prediction $y(s^*)$ at a spatial location $s^*$ is obtained as a weighted average, where observations at closer locations are given higher weights than observations at distant locations. The parameters $\theta(s)$ are defined as a latent stationary Gaussian field. This is equivalent to assuming that $\theta$ has a multivariate normal distribution with mean $\mu$ and spatially structured covariance matrix $\Sigma$, whose generic element is $\Sigma_{ij} = \mathbb{C}\mbox{ov}(\theta(s_i), \theta(s_j))$ \citep{guttorp2006studies, banerjee2003hierarchical,
cressie2015statistics,rue2005gaussian}. 

The concept of spatial process can be extended to the spatio-temporal case including a time dimension. The data are then defined by a process
$$
y(s,t) = \{y(s,t), (s,t) \in \mathcal{D} \times \mathcal{T} \subset \mathbb{R}^d \times \mathbb{R}\}
$$
and are observed at $n$ spatial locations or areas and at $T$ time points. When spatio-temporal geostatistical data is considered \citep{gelfand2010handbook,cressie2015statiotempo}, the kriging model can be extended by defining a spatio-temporal covariance function given as $\mathbb{C}\mbox{ov}(\theta(s_i,t_{\ell}), \theta(s_j,t_h))$. Assuming stationarity in space and time, the space-time covariance function can be written as a function of the spatial Euclidean distance $d(s_i,s_j)$ and of the temporal lag $|t_{\ell}-t_{h}|$. 

In practice, to overcome the computational complexity of these models, some simplifications are introduced. For example, under the separability hypothesis the space-time covariance function is decomposed into the sum (or the product) of a purely spatial and a purely temporal term, e.g. $\mathbb{C}\mbox{ov}(\theta(s_i, t_{\ell}), \theta(s_j, t_h)) = \sigma^2 C_1(d(s_i,s_j)) C_2(|t_{\ell}-t_h|)$ \citep{finkenstadt2006statistical}. Alternatively, it is possible to assume that the spatial correlation is constant in time, giving rise to a space-time covariance function that is purely spatial when $t_{\ell}=t_h$, and is zero otherwise. In this case, the temporal evolution could be introduced assuming that the spatial process evolves in time following an autoregressive dynamics \citep{harvill2010spatio}. Alternatively, several examples of non-separable space-time covariance functions are reported in \cite{cressie1999classes} and \cite{gneiting2002nonseparable}. 

The assumption of Gaussianity in kriging may be inappropriate or difficult to verify, and overly restrictive, e.g. by oversmoothing the prediction. A first relaxation may be represented by finite mixture models. The use of mixture models in the context of spatio-temporal data is useful because it allows to better describe the heterogeneity of the data, which is essential to optimise the allocation of resources and to perform better prediction without the need to include a complete list of covariates (which may not be available). \cite{fernandez2002modelling} and \cite{neelon2014multivariate} propose models for areal (possibly multivariate) data, \cite{frohwirth2008model} introduce a clustering approach for dynamic regression models, while \cite{hossain2014space} analyze space-time mixtures with a particular focus on model selection. However, finite mixture models assume it is possible to consistently select the number of components, which may be difficult, in particular, in presence of covariates. Therefore, extensions to consider infinite mixture models have also been considered. 

Mixture models can be generalised by considering an infinite number of components, through Dirichlet process mixture models. A random distribution $F$ is distributed according to a Dirichlet process (DP) if its marginal distributions are Dirichlet distributed \citep{ferguson1973bayesian,ferguson1974prior}. When analyzing spatial data, we want the Dirichlet process to depend on a covariate represented by the location $s$. Dirichlet processes are then generalised to dependent Dirichlet processes (DDPs) \citep{maceachern2000dependent,barrientos2012support}, which model a collection of unknown distributions indexed by a covariate, through the definition of weights $\pi=(\pi_1, \pi_2, \ldots)$ and atoms $\theta=(\theta_1, \theta_2,\ldots)$ to vary according to the covariate. Spatial dependence can be introduced in the base distribution $F_{0}$, by keeping the weights fixed: this is usually called single-$\pi$ DDP \citep{gelfand2005bayesian,dunson2006semiparametric}. It can also be introduced in the DP weights by keeping the atoms fixed: this is usually called single-atom DDP \citep{reich2007multivariate,dunson2008kernel}. Finally it can be introduced in both the atoms and the weights \citep{duan2007generalized}. 

Introducing the possibility of analysing data evolving over time and space incurs in additional modelling and computational challenges. For this reason, many proposed approaches assume that the DDP only depends on time or space. This work presents a novel Bayesian nonparametric method applied towards problems of inference and prediction in spatio-temporal settings. The method relies on the stick-breaking representation of DPs, with a kernel-based method along the line of \cite{dunson2008kernel}. The model takes advantage of the natural clustering property of the DP to construct clusters evolving with time and space for spatio-temporal data and is able to describe and use for prediction the heterogeneity of the observed data, without the need to introduce a large number of covariates. Observations are allocated so that data points which are observed at nearby locations and at consecutive times are assumed to be more similar. In practice, this means that the model is able to adapt to rapid changes over time and space, to learn more from points close in space and time to improve the prediction accuracy. The case where the DP atoms are not varying with the location and the time allows to define computationally simple and efficient algorithms; extensions to consider both varying weights and varying atoms are conceptually straightforward, but come at an increased computational cost. Both these versions will be considered in this work.

The remaining of the paper is organized as follows: Section \ref{sec:spatialSB} reviews the DP and the stick-breaking process, focusing on the introduction of spatial and temporal dependence; the spatio-temporal proposal, which is the main contribution of this work, is introduced in Section \ref{sec:spatiotempoSB}. The performance of the model is tested on simulation studies in Section \ref{sec:simu} and on a real-data example in Section \ref{sec:real}. Section \ref{sec:conclu} concludes the paper. 

\section{Spatial and Temporal Dirichlet processes}
\label{sec:spatialSB}

Dirichlet processes have been widely applied in Bayesian nonparametric statistics due to their property of conjugacy, their flexibility, and ease of implementation. Taking $F_0$ to be a distribution on $\Theta$, and some real positive number $\alpha$, consider a finite number of sets $A_1, \ldots, A_r \in \Theta$; then $F$ is considered to have a DP prior with base distribution $F_0$ and concentration parameter $\alpha$, i.e. $F \sim DP(\alpha,F_0)$, if
$$
\left( F(A_1), \ldots, F(A_r) \right) \sim Dir(\alpha F_0(A_1), \ldots, \alpha F_0(A_r)).
$$ 
As $\alpha \rightarrow \infty$, $F(A) \rightarrow F_0(A)$, for $A \subset \Theta$. The properties of the DP are clear from its stick-breaking construction \citep{sethuraman1994constructive}, for which, given a concentration parameter $\alpha>0$ and a non-atomic base distribution $F_0$ on $(\Theta, \mathcal{B}(\Theta))$, where $\Theta$ is the parameter space and $\mathcal{B}(\Theta)$ the corresponding $\sigma$-algebra, we have that  $F \sim DP(\alpha, F_0)$ if
\begin{align*}
F(\cdot) &= \sum_{k=1}^{\infty} \pi_k \delta_{\theta_k}(\cdot), \\
\pi_1 &= V_1 \qquad \mbox{and} \qquad
\pi_k = V_k \prod_{j=1}^{k-1} (1-V_j) \quad \mbox{for } k=2, \ldots \\
V_k &\sim Beta(1, \alpha) \qquad \perp \qquad \theta_k \sim F_0 \qquad k=1,2, \ldots \\
\end{align*}
where $\delta_{z}$ is a probability measure concentrated at $z$ (a Dirac mass). This construction is not limited to the Dirichlet process, but generalisation of the beta distribution, such as $V_k \sim Beta(a_{k}, b_{k})$, are also available \citep{ishwaran2001gibbs}. 

In problems of density estimation, in many situations data may be better described by a continuous distribution, therefore the mixing proposed by \cite{antoniak1974mixtures} is usually applied, to produce an infinite mixture of (often Gaussian) kernels $P(\cdot)$, so that $y_i \sim \sum_{k=1}^{\infty} \pi_k p(y_i | \theta_k)$ for $i=1, \ldots, n$, i.e. the model for $y_i$ is assumed to be an infinite mixture model, where the weights of the mixture are defined through the stick-breaking construction. 

Consider a random field $\{y(s): s \in \mathcal{D} \}$ for $\mathcal{D} \subseteq \mathbb{R}^2$, and let $s=(s_1, \ldots, s_n)$ be a set of distinct locations in $\mathcal{D}$. Let $y(s)$, the observation at location $s$, be characterized by the model in \eqref{eq:spatmod}.
The spatial effect $\theta(s)$ is assumed to follow an unknown distribution $F_s$, such that $F_{\mathcal{D}} = \{ F_s: s \in \mathcal{D}\}$ and $F_{\mathcal{D}} \sim \mathcal{P}$ where $\mathcal{P}$ is a probability measure on $(\Gamma,\mathcal{C})$, $\Gamma$ is the space of uncountable collections of probability measures on the Polish space $(\Theta,\mathcal{B})$ indexed by $s \in \mathcal{D}$ and $\mathcal{C}$ is the corresponding $\sigma$-algebra. Therefore $F_s$ is a probability measure over a measurable Polish space $(\Theta, \mathcal{B})$. Then, the density of $y(s)$ is 
$$
g(y(s)|x(s)) = \int_{\Theta} p\left(y(s) | \theta(s) + x(s)^T\beta\right) dF_s.
$$
where $p(\cdot)$ is some kernel, possibly Gaussian. 

The spatial generalisation of the DP can be introduced in the definition of the base distribution, such that $F_{0\mathcal{D}}$ is indexed by $s$, and the atoms $\theta_k$ of the DP are replaced with $\theta_{k\mathcal{D}} = \{ \theta_k(s): s \in \mathcal{D}\}$, for $k=1, 2, \ldots$. This approach can be considered as a particular case of the approach based on regression in the base measure first proposed by \cite{cifarelli1978nonparametric}. For instance, the spatial Dirichlet process (SDP) introduced by \cite{gelfand2005bayesian} uses a base measure $F_{0\mathcal{D}}$ which is a Gaussian process with a (possibly isotropic) covariance function.  
Then $F_{\mathcal{D}}$ is defined such that:
$$
F_s(\cdot) = \sum_{k=1}^{\infty} \pi_k \delta_{\theta_k(s)}(\cdot)
$$ 
for each $s\in \mathcal{D}$. The distribution of the observations $y(s) = (y(s_1), \ldots, y(s_n))$ is almost surely a mixture of Gaussian distributions. This is a model in the class of single-$\pi$ DDPs \citep{maceachern2000dependent}, and maintains the same weights $(\pi_{1}, \pi_2, \ldots, )$ for all locations $s \in \mathcal{D}$. This means that it is difficult to identify the spatial contribution of each component of the mixture. Examples of this approach can be found in  \cite{kottas2007bayesian}, \cite{kottas2012spatial}, and \cite{park2022spatio}.

One of the drawbacks of this approach is that, since the locations associated with the sampled sites constitute one single observation from the random field, replications are needed in order to conduct inference. Such replications are usually taken to be temporal replications. A dynamic model can be introduced in the DP atoms $\theta_{k\mathcal{D}}$ \citep{kottas2008modeling}, however this represents an extension of the modelling of the DP parameters, and the clustering probabilities are still not let evolving with time. \cite{rodriguez2008bayesian} considered a related model, where the atoms are allowed to evolve in time according to an autoregressive model. 

When observations are indexed by time instead of space, \cite{dunson2006bayesian} propose a dynamic mixture of Dirichlet processes (DMDP), where a latent variable dynamically varies across groups; \cite{zhu2005time}, \cite{caron2007bayesian} and \cite{caron2017generalized} propose a time-varying DPM with dynamic evolution of the atoms; \cite{rodriguez2008bayesian} propose a model where atoms are allowed to evolve over time. \cite{lau2008bayesian} consider a model for temporal observations, where each atom can be expressed as an infinite mixture of autoregressions of order $p$. \cite{di2013simple} propose a single-$\pi$ DDP where the atom processes are expressed as linear autoregression $\theta_k(t) = \beta_k + \alpha_k \theta_k(t-1)$. They also consider the case where the atom process is defined as a Ornstein-Uhlenbeck process. An interesting variation of a dynamic DDP construction is proposed by \cite{ascolani2021predictive} who define a family $\mathcal{F} = \{F_t, t \geq 0\}$ of dependent random probability measures $F_t$ which share some, but not all atoms and the set of atoms which are shared in $F_t$ is defined as a pure death process over time, in particular a Fleming-Viot process.  

Models that incorporate dependent weights have a number of theoretical and practical advantages over models that only use dependent atoms. For example, models with non-constant weights have richer support and it has been noticed that constant-weights models cannot generate a set of independent measures \citep{maceachern2000dependent}. Therefore, a second characterization of the DDP for spatial problems relies on the definition of a spatial stick-breaking prior, where the stick-breaking weights are allowed to vary with the spatial locations. The spatial stick-breaking model was first proposed by \cite{reich2007multivariate} and, more generally, by \cite{dunson2008kernel}. We will extend this approach, to address prediction problems including temporal dependence. 
The spatial stick-breaking prior assumes $y(s)$ to be a set of observations at site $s$ such that model \eqref{eq:spatmod} is still valid. 
A prior is assigned to the spatial random effects $\theta(s) \sim F_s$. Differently from the SDP of \cite{gelfand2005bayesian}, this framework utilizes marginal models $F_s, F_{s'}$ rather than a joint model $F_{s,s'}$. In fact, the distribution for $\theta(s)$ is given by 
\begin{align*}
F_s(\cdot) &= \sum_{k=1}^{\infty} \pi_k(s) \delta_{\theta_k}(\cdot) \qquad s \in \mathcal{D} \qquad \mbox{where} \\
\pi_1(s)&=V_1(s) \qquad \mbox{and} \qquad \pi_k(s) = V_k(s) \prod_{j=1}^{k-1} (1-V_j(s))\qquad \mbox{for } k=2,\ldots \\
V_k(s) &= w_k(s,\psi) V_k \\
V_k &\sim Beta(a,b) \qquad \perp \qquad \theta_k \sim F_0 \qquad \mbox{for } k=1,2,\ldots
\end{align*}
The distributions $F_s$ and $F_{s'}$ are related through the dependence among $\pi_k(s)$ and $\pi_k(s')$ which arises from dependence among $V_k(s)$ and $V_k(s')$, according to the bounded functions $w_k(s,\psi)$, which are bounded kernels $w: \mathcal{D} \times \mathbb{R}^2 \rightarrow [0,1]$, allowing the masses of the DP to vary with space. 
An example of kernel function is the squared exponential kernel
\begin{equation}
w_k(s, \psi) = \prod_{j=1}^2 \exp\left[ - \frac{(s_j - \psi_{jk})^2}{h_{jk}^2}\right].
\label{eq:expokern}
\end{equation}
The kernel function is centred around knots $\psi_k = (\psi_{1k},\psi_{2k})$ for $k=1, 2, \ldots$, such that $\psi \sim Q$, where $Q$ is a probability measure on the Polish space $(\Psi, \mathcal{A})$, where $\Psi$ is a Lebesgue-measurable subset of $\mathbb{R}^2$ (not necessarily corresponding to $\mathcal{D}$) and $\mathcal{A}$ is the corresponding Borel $\sigma$-algebra.
The kernel spread is controlled by a bandwidth parameter $h_k = (h_{1k}, h_{2k})$ which can be fixed to a value $\nu$ or allowed to vary according to a prior distribution, e.g. an inverse gamma distribution. The kernel function is defined so that knots closer to $s$ and having a smaller index $k$ have, on average, higher probability weight. $F_s$ and $F_{s'}$ allocate similar probabilities to atoms $\theta=(\theta_1, \theta_2,\ldots)$ if $s$ and $s'$ are close in space. An application of the spatial stick-breaking prior was discussed by \cite{reich2011sufficient} and \cite{reich2011spatial}, within the field of landscape genetics. \cite{reich2007multivariate} and \cite{dunson2008kernel} show that the covariance of $(y(s), y(s'))$ is stationary when integrating with respect to $(V_k, \psi_k, h_k)$, however the conditional covariance can be non-stationary, so this model can be more robust than standard kriging. Moreover, the correlation function tends to one as $s \rightarrow s'$. 

When the focus is on a distribution indexed on time, $F_t$, \cite{mena2016dynamic} construct a single-atom dynamic DDP by setting up a Wrights-Fisher diffusion on the fractions $V_{k}(t)$ in the stick-breaking construction of the marginal DP prior for $F_t$. \cite{gutierrez2016time} propose a time-dependent Bayesian nonparametric model, where the stick-breaking weights are defined on beta variables evolving over time. In discrete time, \cite{dunson2006bayesian} and \cite{griffin2009time} propose an AR-process type model allowing for efficient inference, such that the marginal distribution of $V_{k}(t)$ is $Be(a(t),b(t))$ for all $k$ and $t$. An alternative construction relied on the arrivals contribution of the $m$-DPP \citep{griffin2006order}.

There are few DDP models in the literature where both weights and atoms are allowed to vary according to a covariate. \cite{duan2007generalized} describe a model where atoms $\theta_k(s)$ are i.i.d. from $F_0$ and the weights $\pi_{k(s_1), \ldots, k(s_n)}$ determine the site-specific joint probabilities, however the contribution of time and space to the specific component through the mixing weight is integrated out.  \cite{rodriguez2011nonparametric} propose an approach based on a probit representation of the weights
\begin{align*}
F_s(\cdot) &= \sum_{k=1}^K \pi_k(s) \delta_{\theta_k(s)}(\cdot) \\ 
\pi_k(s) &= \Phi(\alpha_k(s)) \prod_{j < k} (1-\Phi(\alpha_j(s))),
\end{align*}
where $\alpha_k(s)$ have Gaussian marginals and $ \theta_k(s)$ are independent and identically distributed sample paths from a given stochastic process. \cite{arbel2016bayesian} and \cite{barcella2017comparative} empirically show the good predictive performance of models using a probit representation of the weights. Nevertheless, spatial or spatio-temporal dependence among the parameters $\alpha_k$, and indirectly the mixture weights, can be introduced by assuming a Gaussian process over $\mathcal{D}$, which means that the computational cost increases and these models may be difficult to fit in cases with large sample sizes, as the examples in Section \ref{sec:real}. 

To the best of our knowledge, the only work introducing a spatio-temporal varying kernel is proposed by \cite{hossain2013space}, who, however, use the logistic normal kernel as in \cite{fernandez2002modelling}, where spatio-temporal dependence is introduced through spatio-temporal varying covariates in a logistic function. This approach has limitations, including the need for covariates (which may not be available), and the need for a truncation of the possible number of weights (and atoms). In general, this approach has a more difficult interpretation of the interaction of the spatial and temporal components. The contribution of this work is to introduce spatio-temporal dependence in the clustering probabilities, in order to increase the prediction accuracy in highly heterogeneous phenomena.   

\section{Spatio-temporal stick-breaking process}
\label{sec:spatiotempoSB}

In this work, we suppose the observation $y(s,t)$ is a response variable observed at time $t$, with $t \in \{1,\ldots, T\} = \mathcal{T}$, where $T < \infty$ and $\mathcal{T}$ is s distinct set of points in time, and location $s \in \mathcal{D} \subset \mathbb{R}^2$, where $\mathcal{D}$ is the spatial region of interest. While it is easy to define the stick-breaking weights in the case of continuous time, the complete model would imply important modifications and we leave it for further research; see, for example, \cite{mena2011geometric} and \cite{arndt2019sequential}.

Here we propose a class of stick-breaking processes that can be seen as a prior for $F_{\mathcal{D,T}}$ and induces a spatio-temporal dependence on the DP weights. Similarly to model \eqref{eq:spatmod}, $y(s,t)$ is modelled according to 
\begin{equation}
y(s,t) = \theta(s,t) + x(s,t)^T\beta + \varepsilon(s,t) 
\label{eq:spatt_mod}
\end{equation}
where $\theta(s,t)$ is a spatial effect that can evolve over time, $x(s,t)$ is a column-vector of $p$ covariates which may be dependent on the location $s$ and the time $t$ and $\beta$ is a column-vector of $p$ regression coefficients; finally, $\varepsilon(s,t)$ are i.i.d. normal errors. 

The spatio-temporal effect $\theta(s,t)$ is assumed to follow an unknown distribution $F_{s,t}$, such that $F_{\mathcal{D},\mathcal{T}} = \{ F_{s,t}: s \in \mathcal{D}, t \in \mathcal{T}\}$. $F_{\mathcal{D},\mathcal{T}} \sim \mathcal{P}$ where $\mathcal{P}$ is a probability measure on $(\Gamma,\mathcal{C})$, $\Gamma$ is the space of uncountable collections of probability measures on the Polish space $(\Theta,\mathcal{B})$ indexed by $s \in \mathcal{D}$ and $t \in \mathcal{T}$, and $\mathcal{C}$ is the corresponding $\sigma$-algebra. Therefore $F_{s,t}$ is a probability measure over a measurable Polish space $(\Theta, \mathcal{B})$. In particular
$$
g(y(s,t)|x(s,t)) = \int p\left(y(s,t) | \theta(s,t)+ x(s,t)^T\beta\right) dF_{s,t}
$$
where $p(\cdot)$ is a Gaussian kernel. For ease of notation, from now on, we assume no covariates are available. 

The distribution of the observation $y(s,t)$ can be represented as an infinite mixture. The mixing weights are allowed to vary over time and space. They are built similarly to the spatial stick-breaking: 
\begin{align*}
F_{s,t}(\cdot) &= \sum_{k=1}^{\infty} \pi_{k}(s,t) \delta_{\theta_k}(\cdot) \qquad s \in \mathcal{D}, \, t \in \mathcal{T} \qquad \mbox{where} \\
\pi_{1}(s,t)&=V_{1}(s,t) \qquad \mbox{and} \qquad \pi_{k}(s,t) = V_{k}(s,t) \prod_{j=1}^{k-1} (1-V_{j}(s,t))\qquad \mbox{for } k=2,\ldots \\
V_{k}(s,t) &= w_{k}(s,\psi_k,t, \zeta_k) V_k \\
V_k &\sim Beta(a,b) \qquad \perp \qquad \theta_k \sim F_0 \qquad \mbox{for } k=1, 2,\ldots
\end{align*}
Then, we consider a countable sequence of mutually independent random variables $(\psi_1,\psi_2,\ldots)$ such that $\psi_k \sim Q$ for each $k=(1,2,\ldots)$, where $Q$ is a probability measure on a Polish space $(\Psi, \mathcal{A})$, where $\mathcal{A}$ is a Borel $\sigma$-algebra and $\Psi$ is a Lebesgue measurable subset of $\mathbb{R}^2$ which may or may not correspond to $\mathcal{D}$; moreover, we consider a countable sequence of mutually independent random variables $(\zeta_1,\zeta_2,\ldots)$ such that $\zeta_k \sim L$ for each $k=(1,2,\ldots)$ and where $L$ is a probability measure on a Polish space $(\Xi, \mathcal{E})$, where $\Xi$ is a Lebesgue measurable subset of $\mathbb{R}$ and $\mathcal{E}$ is the corresponding Borel $\sigma$-algebra. The random measures are marginally (i.e. for every possible combination of $(s,t) \in \mathcal{S} \times \mathcal{T}$) DP-distributed random measures. The process has full weak support and the DP mixture model induced by the spatio-temporal stick-breaking process is characterized by smooth trajectories as $(s,t)$ ranges over $\mathcal{S} \times \mathcal{T}$ \citep{barrientos2012support}. 

Similarly to the spatial stick-breaking, the spatio-temporal stick-breaking is proper if $\mathbb{E}[V_k]$ and $\mathbb{E}[w_{k}(s,\psi_k,t,\zeta_k)]$ are both positive. These are not restrictive requirements: $V_k$ is usually assumed beta distributed and the prior distribution for $a$, $b$, $\psi_k$ and $\zeta_k$ can be chosen to assure this feature. We propose $V_k \sim Be(a,b)$, where $a$ and $b$ are fixed or can be assumed to follow prior distributions (for example, uniform distributions in a reasonable compact space), the spatial knots $\psi_k=(\psi_{k,1},\psi_{k,2})$ can be fixed, or given uniform prior distributions over the spatial domain, and, similarly, the temporal knots $\zeta_k$ can be fixed or given uniform prior distributions over the temporal domain. 

$F_{s,t}$ and $F_{s',t'}$ allocate similar probabilities to elements of the set $\theta_k$ if $s \rightarrow s'$ and $t \rightarrow t'$, and in this way the process induces dependence on the response variable. $F_{s,t}$ is well defined with $\sum_{k=1}^{\infty} \pi_{k}(s,t) = 1$ almost surely for each $(s,t) \in \mathcal{S} \times \mathcal{T}$ \citep{ishwaran2003generalized}. In order to deal with the sparseness that may be present in the data, it is possible to choose the hyper-parameters $(a, b)$ so that they favour values of $V_k$ close to one or a value (or prior distribution) for the DP concentration parameter $\alpha$ so that the DP assigns high mass to few atoms. 

When dealing with spatial and spatio-temporal data, it is important to study the spatio-temporal correlation or covariance functions assumed by the model. First, we derive the covariance function conditionally on the beta variables $V_k$, the knots $\psi_k$ and $\zeta_k$. Since the spatio-temporal stick-breaking process has discrete realizations, the covariance function is:
\begin{small}
\begin{align*}
\mathbb{C}&\mbox{ov}(y(s,t), y(s',t')|V, \psi, \zeta) = \Pr(\theta(s,t) = \theta(s',t') | V, \psi, \zeta) = \sum_{k=1}^{\infty} \pi_{k}(s,t) \pi_{k}(s',t') \\
&=\sum_{k=1}^{\infty} \left[ w_{k}(s,\psi_k,t,\zeta_k) w_{k}(s',\psi_k,t',\zeta_k) V_k^2 \prod_{j < k} [1-w_{j}(s,\psi_j,t,\zeta_j) V_j] [1-w_{j}(s',\psi_j,t',\zeta_j) V_j] \right] \\
&= \sum_{k=1}^{\infty} \left[ w_{k}(s,\psi_k,t,\zeta_k) w_{k}(s',\psi_k,t',\zeta_k) V_k^2  \prod_{j < k} [1- (w_{j}(s,\psi_j,t,\zeta_j)+w_{j}(s',\psi_j,t',\zeta_j)) V_j + \right. \\
& \qquad +\left. w_{j}(s,\psi_j,t,\zeta_j)w_{j}(s',\psi_j,t',\zeta_j) V_j^2]  \right].  
\end{align*}
\end{small}
The covariance function, as expected, depends on the choice of the kernel, and this means that it is maximised when choosing a kernel for which $w(s,\psi,t,\zeta) \rightarrow w(s',\psi,t',\zeta)$ when $s \rightarrow s'$ and $t \rightarrow t'$. It is also interesting to notice that, the conditional covariance function does not depend on the specific value of the response variable. This covariance can be a function of only the distances between spatial locations and between time points or can depend on the specific spatial locations and time points, therefore the fact that the process is stationary depends on the choice of kernel. 

To better understand the properties of the covariance function between $y(s,t)$ and $y(s',t')$, it is useful to marginalise over the random beta variables $V_k$, the knots and the kernel parameters to obtain the unconditional covariance. It is reasonable to use independent priors: for example, uniform independent priors for $\psi_k$ over the spatial domain and uniform independent priors for $\zeta_k$ over the temporal window of the observations. Then, the unconditional covariance function is 
\begin{small}
\begin{align*}
\mathbb{C}\mbox{ov}(y&(s,t), y(s',t')) = \mathbb{E}_{[V, \psi,\zeta]} \left\{ \sum_{k=1}^{\infty} \left[ w_{k}(s,\psi_k,t,\zeta_k) w_{k}(s',\psi_k,t',\zeta_k) V_k^2 \cdot \right. \right. \\
& \left. \left. \prod_{j < k} [1- (w_{j}(s,\psi_j,t,\zeta_j)+w_{j}(s',\psi_j,t',\zeta_j)) V_j + w_{j}(s,\psi_j,t,\zeta_j)w_{j}(s',\psi_j,t',\zeta_j) V_j^2]  \right]  \right\} \\
&= \int_{\Psi} \int_{\Xi} w(s,\psi,t,\zeta) w(s',\psi,t',\zeta) \pi(\psi,\zeta) \mathbb{E}[V_k^2] \cdot \sum_{k=1}^{\infty} \left[ 1- \int_{\Psi} \int_{\Xi} w(s,\psi,t,\zeta) \pi(\psi,\zeta) \mathbb{E}[V_k]  + \right. \\
& \left. - \int_{\Psi} \int_{\Xi} w_{t}(s',\psi,t',\zeta) \pi(\psi,\zeta) \mathbb{E}[V_k] + \int_{\Psi} \int_{\Xi} w(s,\psi,t,\zeta) w(s',\psi,t',\zeta) \pi(\psi,\zeta) \mathbb{E}[V_k^2] \right]^{k-1}.
\end{align*}
\end{small}
If we use $V_k \sim Be(a,b)$ for all $k$, then $\mathbb{E}[V_k]$ and $\mathbb{E}[V_k^2]$ have a fixed form, depending on $a$ and $b$; the unconditional covariance does not depend on the value of $y(s,t)$ or $y(s',t')$ and it is equal to one only when $w(s,\psi,t,\zeta) \rightarrow w(s',\psi,t',\zeta)$ when $s \rightarrow s'$ for every $\psi \in \Psi$ and $t \rightarrow t'$ for every $\zeta \in \Xi$. Moreover, the covariance does not depend on $F_0$. The covariance structure depends on the expectations of the kernels, i.e. the expectation of $w_{k}(s,\psi_k,t,\zeta_k)$, the expectation of $w_{k}(s',\psi_k,t',\zeta_k)$, and the expectation of $w_{k}(s,\psi_k,t,\zeta_k)w_{k}(s',\psi_k,t',\zeta_k)$, therefore it depends on the form of the chosen kernel. 
It is important to notice that, similarly to all single-atom DDPs, the covariance function presents a lower bound. However, this does not influence the fact that the process has full support with respect to the weak topology, and the corresponding DP mixture model has large Hellinger support \citep{barrientos2012support}.

\begin{figure}
\centering
\begin{subfigure}{.45\textwidth}
  \centering
  \includegraphics[width=5.5cm,height=5cm]{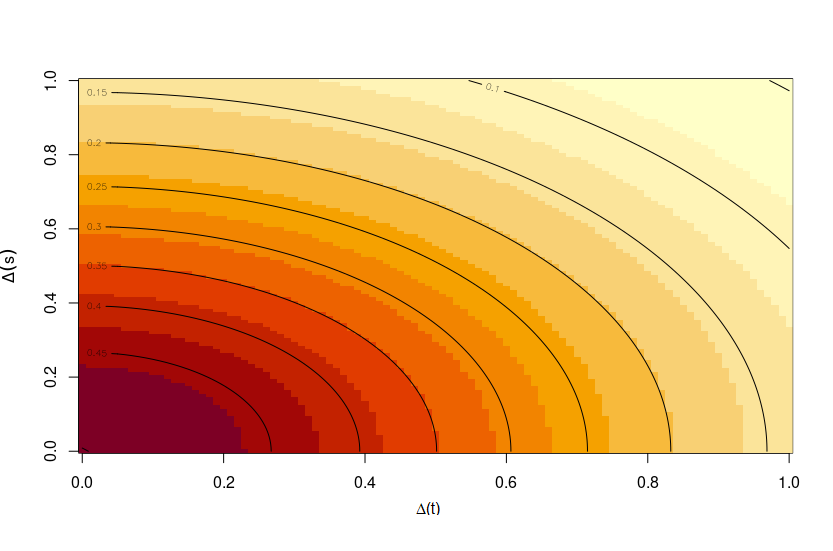}
\end{subfigure}%
\begin{subfigure}{.45\textwidth}
  \centering
  \includegraphics[width=5.5cm,height=5cm]{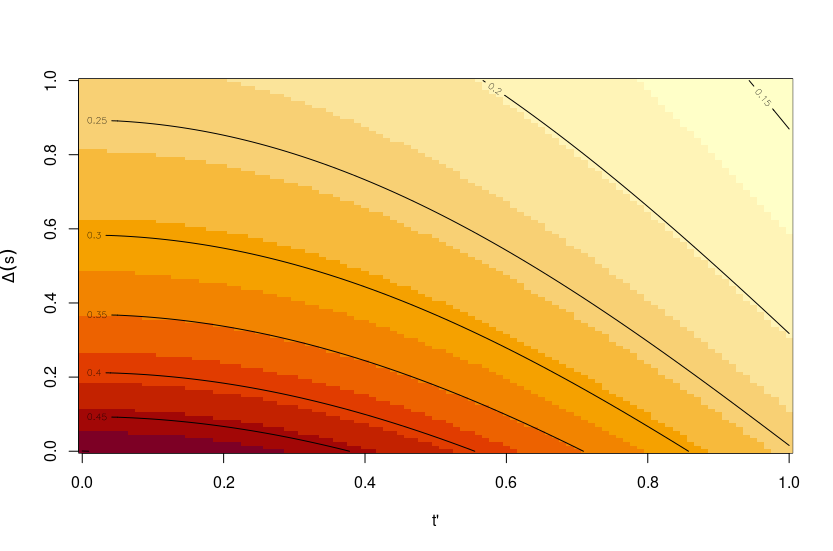}
\end{subfigure}
\caption{Covariance function of the spatio-temporal process when separable exponential kernels are used for space and time (left), and when the non-separable kernel \eqref{eq:gneiting} is used (right). In the second case, the covariance structure does not only depend on the distance between time points, so in this plot we fix $t=0$, and compute the covariance by varying $t'$.}
\label{fig:corr_sepkernels}
\end{figure}

The kernel $w(s,\psi,t,\zeta)$ is a bounded kernel, $w: \mathcal{D} \times \mathbb{R}^2 \times \mathcal{T} \times \mathbb{R} \rightarrow [0,1]$, allowing masses to vary over time and space. The easiest way to define a spatio-temporal kernel is by considering a separable kernel
$$
w_k(s, \psi_k, t,\zeta_k) = w_k(s, \psi_k) \cdot w_k(t, \zeta_k) \qquad k=1,2,\ldots 
$$
The hypothesis of separability in this setting implies that the temporal autocorrelation of the mixing weights for two locations at distance $d$ is proportional to the temporal autocorrelation of the mixing weights for two locations at distance $rd$ for some constant $r$. While the kernel is separable, the covariance structure of the process is not separable. Indeed, the covariance function of the process is defined as
\begin{small}
\begin{align}
\mathbb{C}\mbox{ov}&(y(s,t), y(s',t') | \psi, \zeta) = \Pr[\theta(s,t) = \theta(s',t')|\psi, \zeta] = \nonumber \\
& \sum_{k=1}^{\infty} \left[ w_k(s, \psi_k) w_k(t, \zeta_k) w_k(s', \psi_k) w_k(t', \zeta_k) V_k^2 \right. \nonumber \\
& \qquad \left. \cdot \prod_{j<k} \left( 1-( w_j(s,\psi_j) w_j(t, \zeta_j) + w_j(s', \psi_j) w_j(t', \zeta_j) ) V_j  + \right. \right. \nonumber \\
& \qquad \qquad \left. \left. w_j(s, \psi_j) w_j(t, \zeta_j) + w_j(s', \psi_j) w_j(t', \zeta_j)  \right) \right] \label{eq:covsepkern}
\end{align}
\end{small}
and, by marginalising with respect to $\psi$ and $\zeta$,
$$
\Pr(\theta(s,t) = \theta(s',t')) = \frac{g(s,s')g(t,t')}{2\left(1+\frac{b}{a+1} - g(s,s')g(t,t')\right)}
$$
where $a$ and $b$ are the parameters of the beta distribution such that $V_k \sim Be(a,b)$ and 
\begin{align*}
g(s,s',t,t') &= \frac{\int_{\Psi} w(s, \psi)w(s', \psi)p(\psi) d\psi}{\int_{\Psi} w(s,\psi)p(\psi) d\psi} \cdot \frac{\int_{\Xi} w(t, \zeta)w(t', \zeta)p(\zeta) d\zeta }{\int_{\Xi} w(t,\zeta)p(\zeta) d\zeta}. 
\end{align*}
The covariance function is not separable, in the sense that it is not a product of functions in time and space, however it is a combination of separable factors. While the model is already flexible, leading to a non-separable covariance function, it implies a symmetry between time and space. For example, consider using an exponential kernel for both time and space as in \eqref{eq:expokern}; Figure \ref{fig:corr_sepkernels} shows that there is a symmetry between time and space, such that the covariance between the observations decreases symmetrically as the distance between time and spaces increases, and this may represent a strong assumption depending on the phenomenon, in particular when this interaction is not known. Importantly, the conditional covariance function is non-stationary, but the unconditional covariance function only depends on the distance between locations and time points, i.e. it is stationary, when choosing separable exponential kernels. 

The impact of the hypothesis of separability among spatio-temporal processes is already unclear; a more complicated problem is to study its impact on the mixing weights, given that it relates to non-observable variables. There are several non-separable kernels proposed in the literature: see, for example, \cite{cressie1999classes}, \cite{gneiting2002nonseparable} and \cite{genton2007separable}. We propose to use the class of kernels proposed by \cite{gneiting2002nonseparable}:
$$
w(s,\psi,t,\zeta) = \frac{1}{\xi(|t-\zeta|^2)^{d/2}} \phi\left( \frac{(s_1 - \psi_1)^2+(s_2 - \psi_2)^2}{\xi(|t-\zeta|^2)}\right)
$$
where $\phi(\cdot)$ is a completely monotone function, and $\xi(\cdot)$ is a positive function with a completely monotone derivative. For example, 
\begin{equation}
w(s,\psi,t,\zeta) = \frac{1}{\gamma |t-\zeta| + 1} \exp\left( - \frac{(s_1 - \psi_1)^2+(s_2 - \psi_2)^2}{(\gamma |t-\zeta| + 1)^{\lambda/2}}\right),
\label{eq:gneiting}
\end{equation}
where $\gamma$ is a non-negative scaling parameter and $\lambda\in [0,1]$ is a parameter controlling the space-time interaction: if $\lambda=0$, the model is separable; as $\lambda \rightarrow 1$, the space-time interaction increases and the spatial dependence at positive lags decreases more slowly. The parameter $\gamma$ can be fixed or given an inverse-gamma or a uniform prior distribution. 

The definition of $\lambda$ allows an automatic testing about it being equal to zero (i.e. separability of space and time in the definition of the mixing weights), by assuming a spike-and-slab prior \citep{mitchell1988bayesian,ishwaran2005spike,andersen2017bayesian}, so that
$$
\lambda \sim \omega_{\lambda} f(\lambda) + (1-\omega_{\lambda}) \delta_{0},
$$
i.e. the prior distribution results in a mixture between a continuous distribution $f(\lambda)$ (for example, uniform or beta) and a Dirac mass at zero $\delta_0$, with some mixing weights $\omega_{\lambda} $ that can be fixed or allowed to follow a prior distribution. In this way, the posterior distribution for $\lambda$ will be again a mixture, and it is possible to test the hypothesis of separability by looking at the weight of the mass component at zero. 

The structure of the covariance function among observations is similar to the one obtained with separable kernels in Equation \eqref{eq:covsepkern}, however the function $g(\cdot, \cdot)$ is not defined on separable factors anymore. More specifically, 
\begin{align*}
g(s,s',t,t') = \frac{\int_{\Psi} \int_{\Xi} w_k(s,\psi_k,t,\zeta_k) w_k(s',\psi_k,t',\zeta_k) p(\psi,\zeta) d\psi_k d\zeta_k}{\int_{\Psi} \int_{\Xi} w_k(s,\psi_k,t,\zeta_k)p(\psi, \zeta) d\psi_k d\zeta_k},
\end{align*}
when considering $\gamma$ and $\lambda$ fixed.  
In the case of kernel \eqref{eq:gneiting}, this function becomes
$$
g(s,s',t,t') = \frac{1}{\gamma|t'|+1} \sqrt{(\gamma |t'| + 1)^{\lambda/2}} \exp\left[ - \frac{(s_1-s_1')^2 + (s_2-s_2')^2}{(\gamma |t| + 1)^{\lambda/2} + (\gamma |t'| + 1)^{\lambda/2}} \right],
$$
and it is clear that this function does not only depend on the distance between time-points. 
Therefore, depending on the choice of the kernel, the covariance can be stationary or non-stationary, and it is non-separable with respect to space and time. In the particular choice of kernel \ref{eq:gneiting}, the covariance function obtained by marginalising over the parameters of the kernel is stationary with respect to space, but not with respect to time. Therefore, the covariance obtained for the spatio-temporal stick-breaking can be non-stationary, differently from the spatial stick-breaking. 

An example of covariance structure of the process is illustrated in Figure \ref{fig:corr_sepkernels}, which shows that the covariance can decrease in a non-symmetric way with respect to the direction of time and space, for instance, by decreasing more slowly along the time direction. Figure \ref{fig:corr_sepkernels} is defined for a fixed value of $\gamma$ and $\lambda$, however, in real situations, the parameters are estimated, therefore the type of interaction between time and space is left to be defined by the data. Supplementary Material shows additional examples of the covariance function of the process by varying the parameters ($\gamma$, $\lambda$, $\psi$, $\zeta$). 

\begin{figure}
\begin{multicols}{2}
    \includegraphics[width=5.5cm,height=4cm]{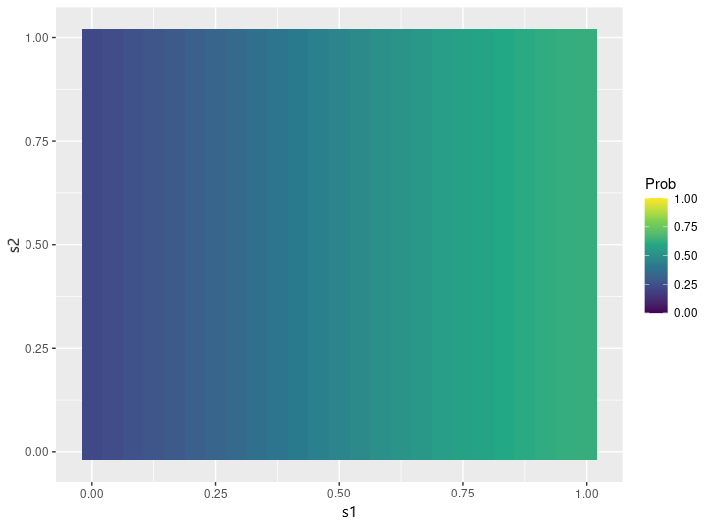}\par 
    \includegraphics[width=5.5cm,height=4cm]{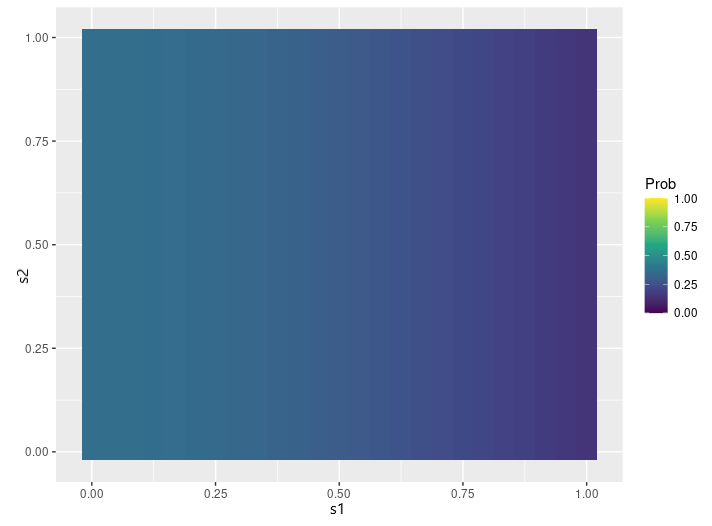}\par 
    \end{multicols}
\caption{Map of $\pi_1$ (left) and $\pi_2$ (right) for each location and $t=1$.}
\label{fig:pi1pi2}
\end{figure}

Observations which are close in space and time are more likely to be allocated to similar components, because the mixing weights for each component are more similar. Figure \ref{fig:pi1pi2} shows a map of the weights $\pi_1(s,t=1)$ and $\pi_2(s,t=1)$ over the space domain, which represents how locations tend to be more likely allocated to component one or two. 

\begin{figure}
\includegraphics[width=7cm,height=5cm]{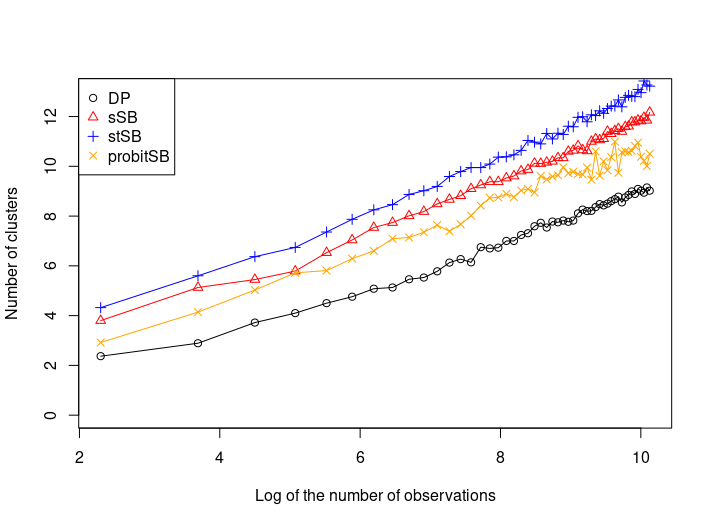}
	\caption{Clustering structure generated by the spatio-temporal stick-breaking. Comparison with the classic DP, the spatial stick-breaking and the probit stick-breaking process. The figure shows the expected number of clusters with respect to the log number of observations.}
\label{fig:numbclust}
\end{figure}

Finally, Figure \ref{fig:numbclust} shows the behavior of the number of allocated clusters of four processes: classic DP, spatial stick-breaking, spatio-temporal stick-breaking and probit stick-breaking. For all the processes, as expected, the number of populated clusters increases with the number of observations. The spatio-temporal stick-breaking on average is associated with a larger number of populated clusters than the other processes: in terms of prediction and density estimation, this property can be associated with a higher ability to take into account local heterogeneity, and can lead to higher prediction scores in more complex situations, as it will be shown in Section \ref{sec:real}. 
  
Under the single-atom setting, the dependence on space and time is expressed through the weights of the stick-breaking representation. One advantage of doing so is that the implied partition model changes with the values of $s \in \mathcal{S}$ and $t \in \mathcal{T}$. This is an essential property, in particular when the allocation of the observations to the mixing component is also of interest. Another important feature is that problems related to inference of $\theta_k(s,t)$ are avoided: this has computational consequences, i.e. the approaches are less computationally intensive, and predictive consequences, i.e. the atoms $\theta_k(s,t)$ should include values which are not linked with any observed data in the case of inference for new values of $(s,t)$. However, extensions are also possible. 

\subsection{Allowing the atoms to depend on space and time}
\label{sub:varyingweights}

The single-atom stick-breaking process is capable of capturing features of distributions that vary over time and space, even whithout the need for replicates for each location in the index space. However, any nonparametric model that does not allow both weights and atoms to vary over the covariate space assumes a lower bound to the covariance between the
distributions \citep{rodriguez2011nonparametric}. This is not necessarily a drawback of the model, if the assumption of dependence is a reasonable one, but it may be an unreasonable assumption in the case we want to take into account the possibility of no dependence over time or space. 

Atoms may also be allowed to vary by introducing a modification of the stick-breaking representation of each element 
\begin{align}
F_{s,t}(\cdot) = \sum_{k=1}^{\infty} \pi_k(s,t) \delta_{\theta_k(s,t)}(\cdot)
\label{eq:wanda}
\end{align}
where $\pi_k(s,t)$ are defined as in the previous section, and $\theta_k(s,t)$ are independent stochastic processes with index set $\mathcal{S} \times \mathcal{T}$ and $F_{s,t,0}$ marginal distributions. The processes associated to the weights and atoms are independent. 

The atoms $\theta_k(s,t)$ can be generated by a Gaussian process, i.e. 
$$
\theta_{\mathcal{S}\times \mathcal{T}} \sim \sum_{k=1}^{\infty} \pi_k(s,t) \delta_{\theta_{k, \mathcal{S}\times\mathcal{T}}},
$$
where $\theta_{k,\mathcal{S}\times \mathcal{T}} = \{\theta_k(s,t), (s,t) \in \mathcal{S}\times\mathcal{T}\}$. The Gaussian process is then characterised by a covariance matrix considering space and time distances. In this work, we use a spatio-temporal covariance function defined as the product of one covariance function depending on space and one covariance function depending on time. More specifically, if $d(s,s')$ is the distance between two spatial locations, for instance the Euclidean distance, and $d(t,t_0)$ is the temporal lag with respect to the initial time point, then the covariance function is 
$$
C(d(s,s'),d(t,t_0)) = C_s(d(s,s')) \times C_t(d(t,t_0)),
$$
where $C_s$ is a spatial covariance function and $C_t$ is a temporal covariance function. In this work, we use an exponential covariance function for $C_s(\cdot)$ and the temporal covariance corresponding to an $AR(1)$ process $C_t(\cdot) = \rho (t-t_0)$. Such covariance matrix is implemented in the \texttt{SpatialTools} package available in \texttt{R} \citep{french2022package}. 

\begin{figure}
\begin{multicols}{4}
	\textbf{$t=1$}\par
    \includegraphics[width=\linewidth]{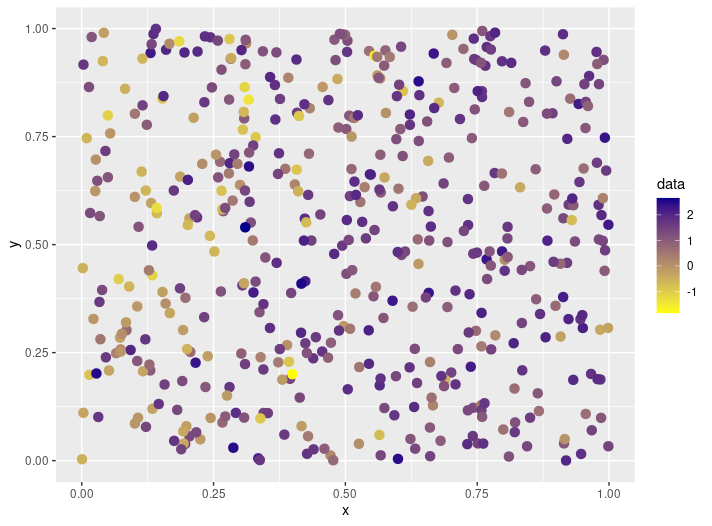}\par 
    \textbf{$t=2$}\par
    \includegraphics[width=\linewidth]{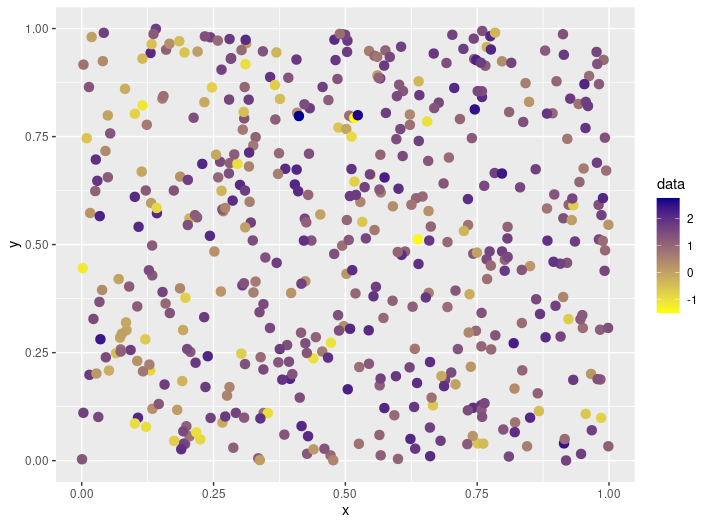}\par 
    \textbf{$t=1$}\par
    \includegraphics[width=\linewidth]{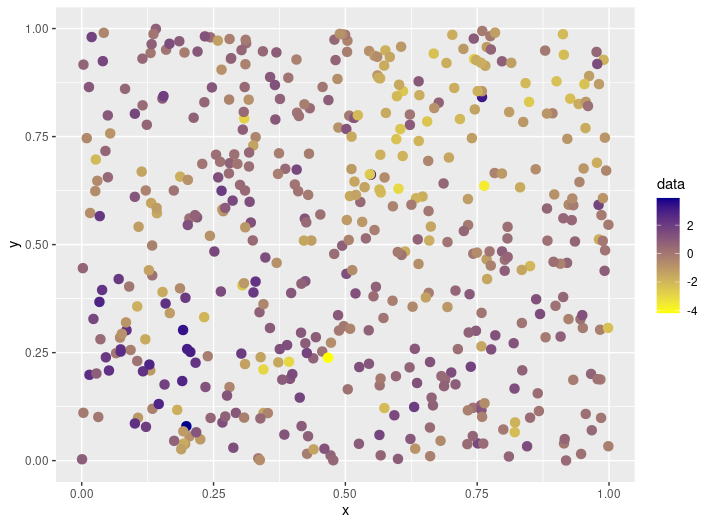}\par 
    \textbf{$t=2$}\par
    \includegraphics[width=\linewidth]{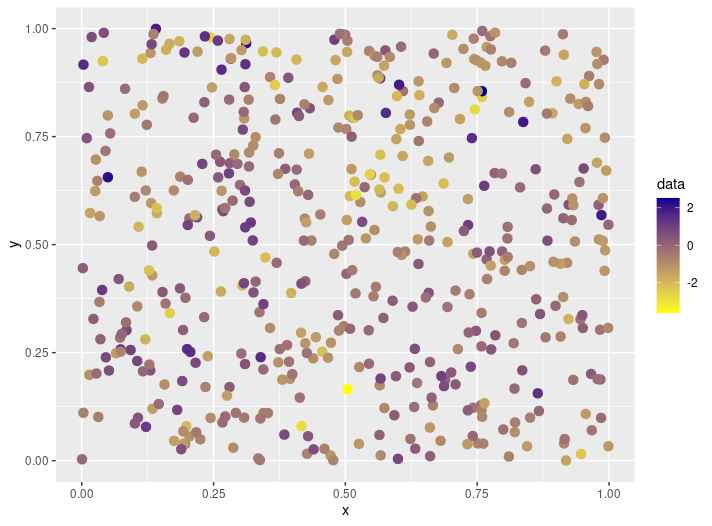}\par 
    \end{multicols}
\begin{multicols}{4}
    \textbf{$t=3$}\par
    \includegraphics[width=\linewidth]{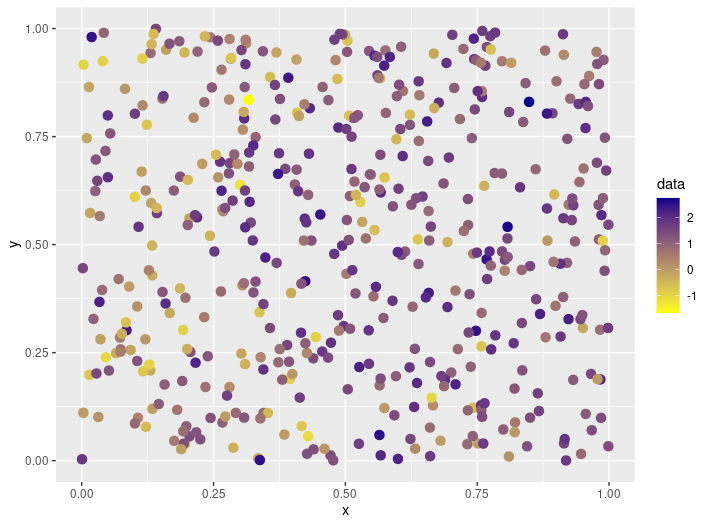}\par 
	\textbf{$t=4$}\par
    \includegraphics[width=\linewidth]{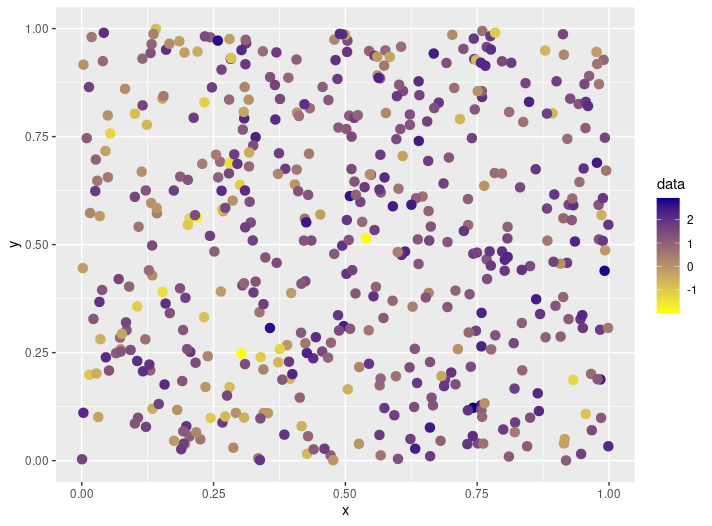}\par 
    \textbf{$t=3$}\par
    \includegraphics[width=\linewidth]{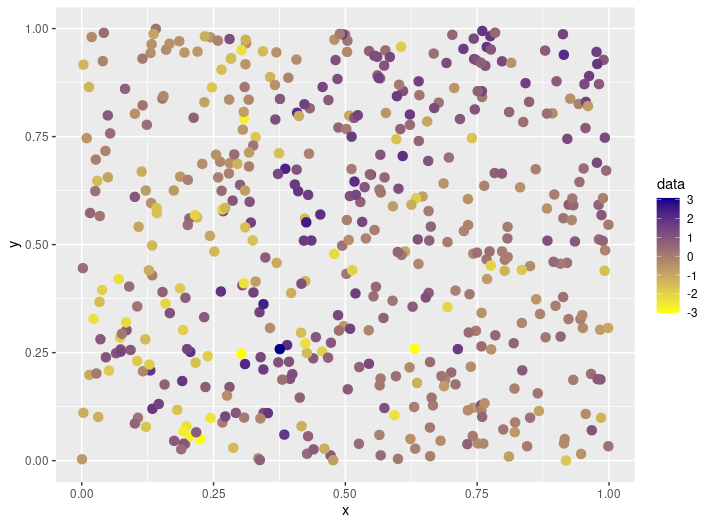}\par 
    \textbf{$t=4$}\par
    \includegraphics[width=\linewidth]{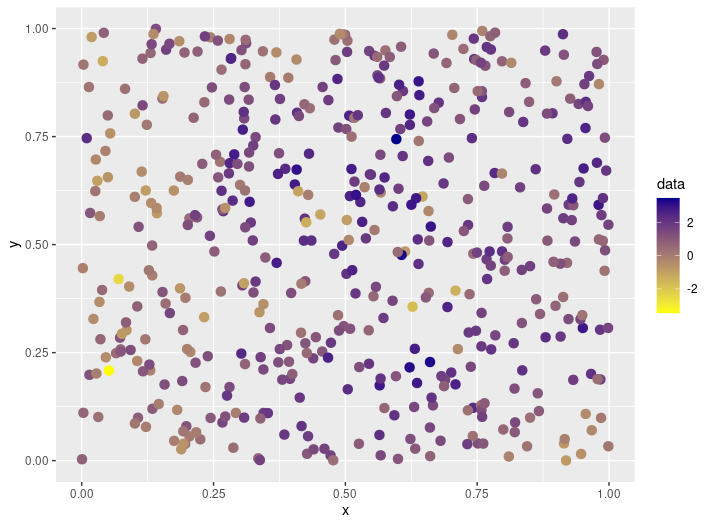}\par 
    \end{multicols}
\begin{multicols}{4}
    \textbf{$t=5$}\par
    \includegraphics[width=\linewidth]{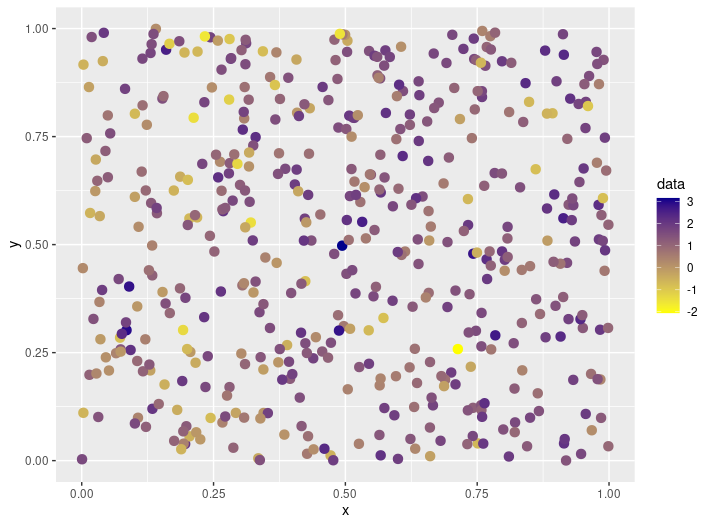}\par 
    \textbf{$t=6$}\par
    \includegraphics[width=\linewidth]{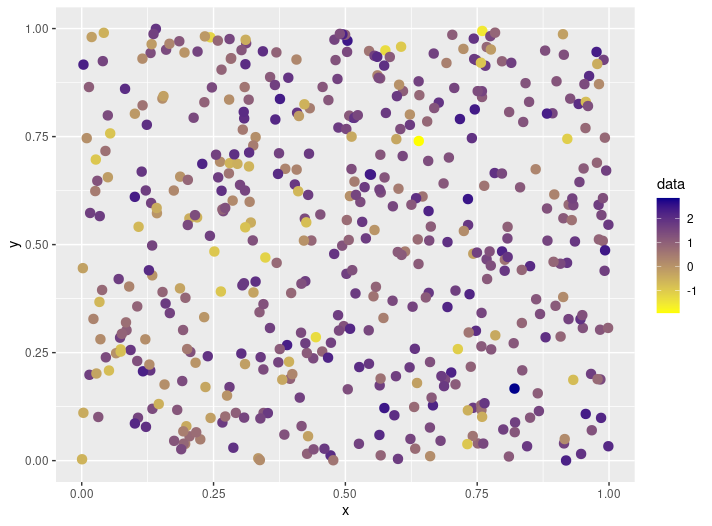}\par 
    \textbf{$t=5$}\par
    \includegraphics[width=\linewidth]{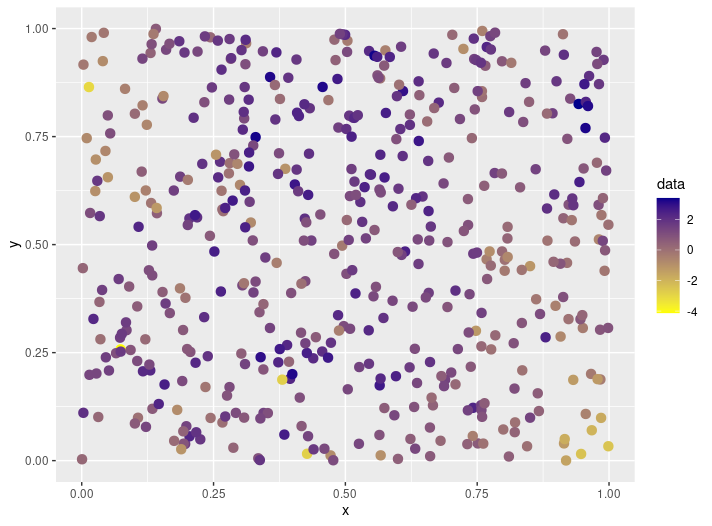}\par 
    \textbf{$t=6$}\par
    \includegraphics[width=\linewidth]{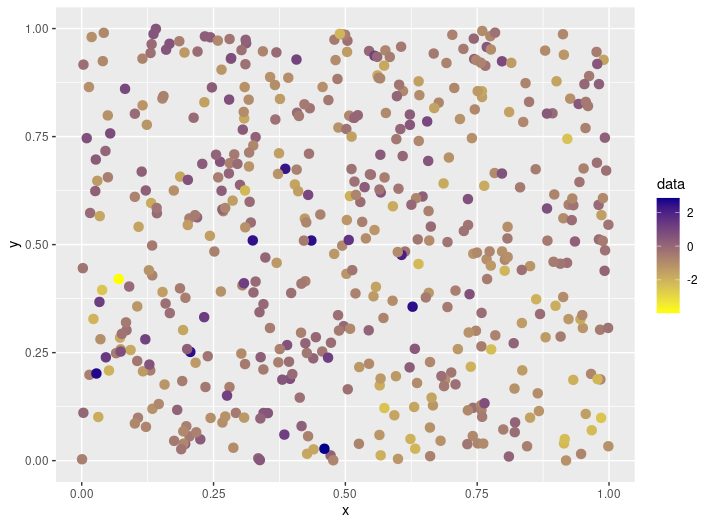}\par 
    \end{multicols}
\begin{multicols}{4}
    \textbf{$t=7$}\par
    \includegraphics[width=\linewidth]{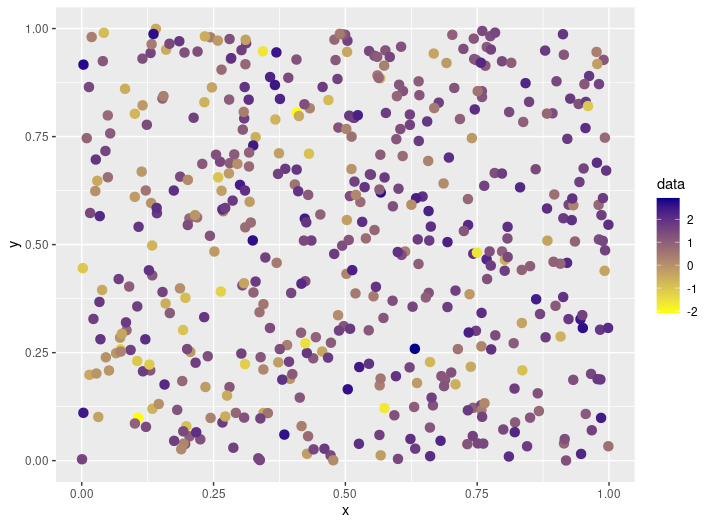}\par 
    \textbf{$t=8$}\par
    \includegraphics[width=\linewidth]{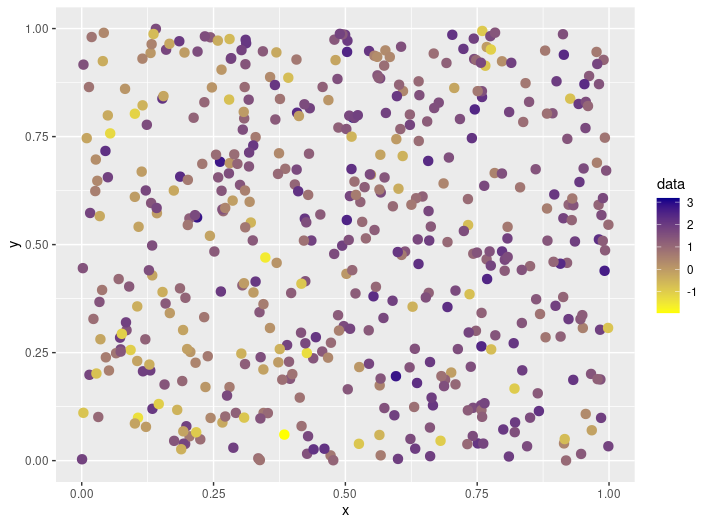}\par 
    \textbf{$t=7$}\par
    \includegraphics[width=\linewidth]{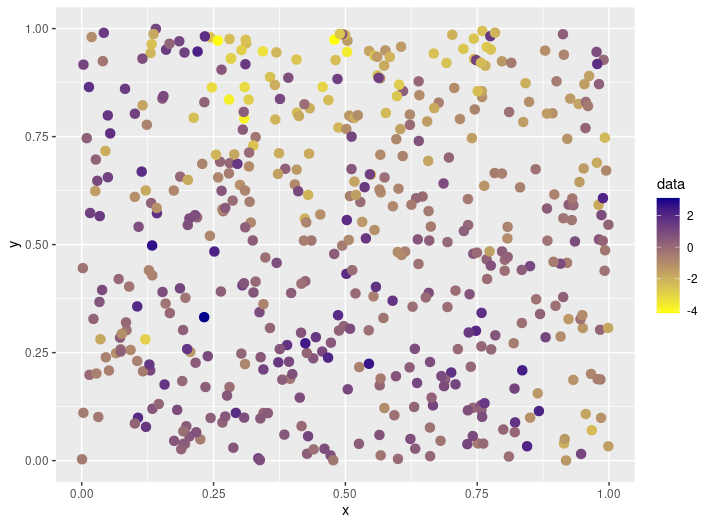}\par 
    \textbf{$t=8$}\par
    \includegraphics[width=\linewidth]{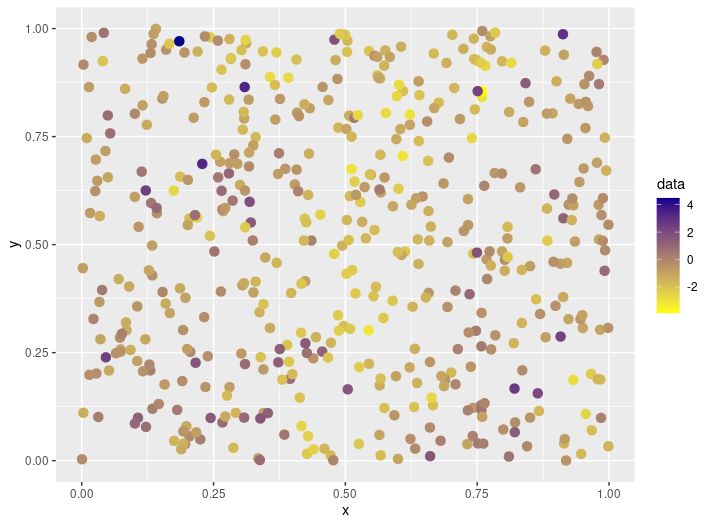}\par 
    \end{multicols}
\begin{multicols}{4}
    \textbf{$t=9$}\par
    \includegraphics[width=\linewidth]{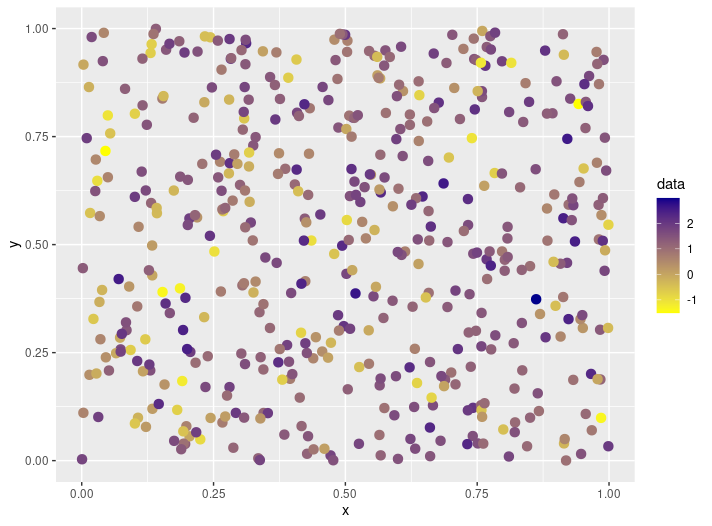}\par 
    \textbf{$t=10$}\par
    \includegraphics[width=\linewidth]{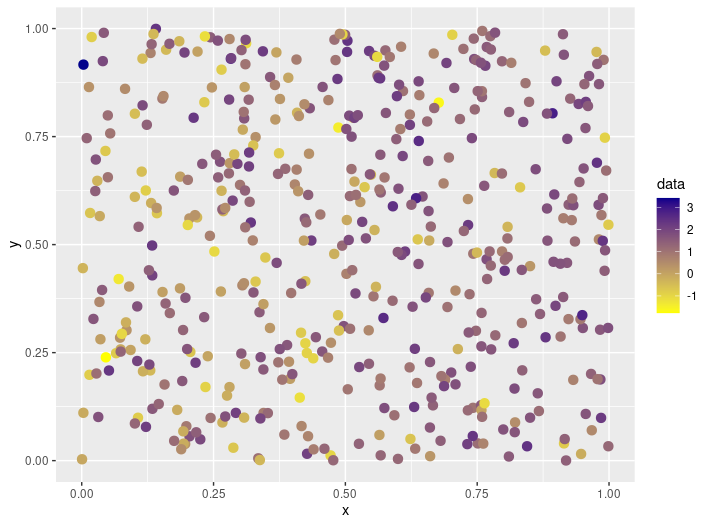}\par
    \textbf{$t=9$}\par
    \includegraphics[width=\linewidth]{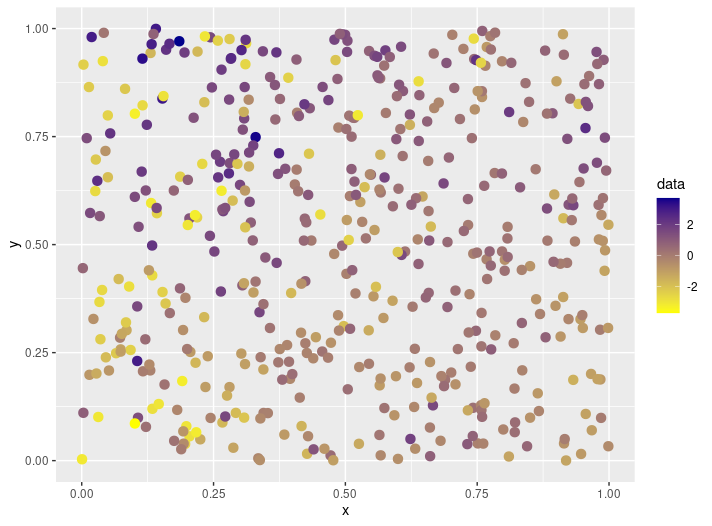}\par 
    \textbf{$t=10$}\par
    \includegraphics[width=\linewidth]{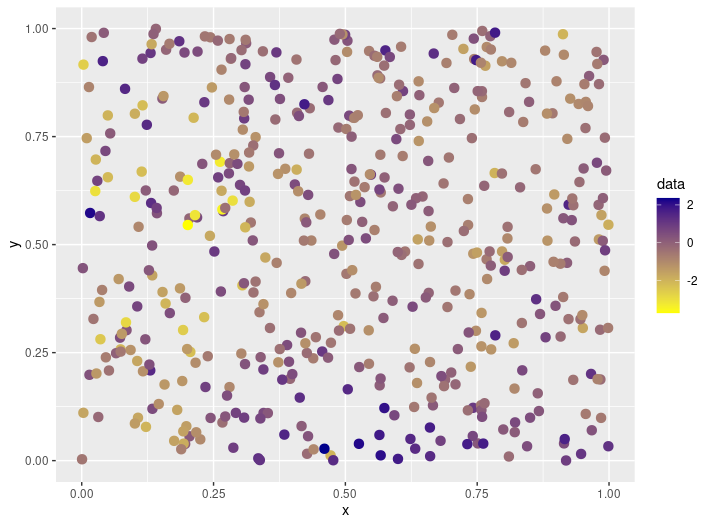}\par 
    \end{multicols}
\caption{Temporal evolution of the observations for uniformly selected locations. Process with only varying weight on the left (columns 1 and 2) and process with both varying weights and atoms on the right (columns 3 and 4).}
\label{fig:simulations1_time}
\end{figure}

Figure \ref{fig:simulations1_time} shows the temporal evolution of the process over time and space; on the left, the single-atom process is represented, and on the right the process with both varying weights and atoms is represented. Figure \ref{fig:simulations1_time} shows a case where locations are uniformly generated in a unit-square set; the Supplementary Material includes additional simulation, considering more structured locations. The process seems to vary more smoothly over space and with a higher variation over time by introducing spatio-temporal dependence also among the atoms with respect to the process only introducing dependence among weights. From the simulations, it seems that the model with both weights and atoms dependent on $(s,t)$ can better describe heterogeneity and situations where the process varies more quickly over time. However, the main drawback of this process is that inference for it is more computationally expensive. For this reasons, in some simulation settings and the real-data applications we fit the model on a subset of the datasets. 

\subsection{Computational details}

Posterior inference can be obtained through Markov Chain Monte Carlo (MCMC). Usually, two different finite approximation strategies are adopted for stick-breaking processes. The first strategy is a marginal approach which uses a P\'olya urn representation to marginalise over the Dirichlet process prior. The second strategy is a conditional approach \citep{ishwaran2000markov,
papaspiliopoulos2008retrospective}, which avoids marginalization. We use this second strategy. We present details of the MCMC scheme for the case of both atoms and weights varying over $(s,t)$, while the single-atom case results as a simplification of this scheme. 

It is possible to use a truncation scheme or to consider the possibility to have an infinite number of components. Earlier works have investigated the truncation level, see, for example, \cite{muliere1998approximating} and \cite{ishwaran2001gibbs}. In our simulation studies and real-data applications, we have set the maximum number of clusters to a very large number, i.e. $M=100$.

In the following, for ease of notation, we consider to observe $\{y_i: i \in (1,2,\ldots,n)\}$ with $y_i = y(s,t)$. Moreover
$$
y_i = \theta_i + \varepsilon_i.
$$
Let $(\theta_1^*,\ldots, \theta_K^*)$ denote $K < n$ unique values for the atoms of the process. Let $c_i = k$ if $\theta_i = \theta_k^*$, i.e. subject $i$ is allocated to component $k$, with $c=(c_1, \ldots, c_n)$, and let $z_k = j$ denotes that atom $\theta^*_k$ comes from $F_j$, with $z=(z_1, \ldots, z_K)$. 

Let $\theta^{-i}$, $\theta^{*-i}$, $c^{-i}$, $z^{-i}$ and $y^{-i}$ correspond to the vectors $\theta$, $\theta^*$, $c$, $z$ and $y$ without the contribution of subject $i$. The number of subjects allocated to component $k$ is $n_k = \sum_{i=1}^n \mathbb{I}[c_i = k]$, with $\sum_{k=1}^{\infty} n_k = n$. Following \cite{pitman1996some} and \cite{ishwaran2003generalized}, 
$$
p(\theta_i \in \cdot | c_i = k,c^{-i},z^{-i},\theta^{*-i}) = w_{ik0} F_0  + \sum_{j \in \mathcal{S}_k^{-i}} w_{ikj} \delta_{\theta_j}(\cdot)
$$
where $\mathcal{S}_k^{-i} = \mathcal{S}_k \cap \{1,2,\ldots, n\}\setminus \{i\}$, $\mathcal{S}_k = \{i : c_i = k, i=1, 2, \ldots, \infty\}$ is the subset of positive integers indexing subjects allocated to component $k$. The weights are defined as $w_{ikj} = \frac{1}{\alpha + |\mathcal{S}_k^{-i}|}$, where $|\cdot|$ stands for the cardinality, and $w_{ik0} = \frac{\alpha}{\alpha + |\mathcal{S}_k^{-i}|}$. The parameter $\alpha$ is the concentration parameter of the DP. By marginalising out $c_i$
$$
\Pr(\theta_i \in \cdot | z^{-i},c^{-i},\theta^{*-i}) = p_{i0} F_0(\cdot) + \sum_{k=1}^{K^{-i}} p_{ik} \delta_{\theta_k^{-i}} (\cdot) + p_{i,b^{-i}+1} F_0(\cdot), 
$$
where $K^{-i}$ is the length of $\theta^{-i}$ and 
\begin{align*}
p_{i0} &= \sum_{k \in \mathcal{I}^{-i}} \pi_{ik} w_{ik0} \\
p_{ik} &= \pi_{ik}^{-i} \sum_{j: c_j^{-i}=k} w_{i z_k^{-i}j} \qquad j=1, \ldots, K^{-i} \\
p_{i,K^{-i}+1} &= \sum_{k \in \mathcal{I}^{C,-i}} \pi_{ik} w_{ik0},  
\end{align*} 
where $\mathcal{I}$ is the set of occupied component, i.e. $\{k \in (1,2,\ldots): n_k > 0\}$, while $\mathcal{I}^C$ is the set of empty components.   

Suppose the likelihood contribution for subject $i$ is indicated by $g(y_i | \theta_i)$, then 
\begin{equation*}
\Pr(c_i = k | y, c^{-i},z^{-i}, \theta^{*-i}) \propto \begin{cases}
\mbox{const}_i \cdot w_{ik} g_0(y_i) \qquad k=0, K^{-i}+1 \\
\mbox{const}_i \cdot w_{ik} g(y_i | \theta_k^{-i}) \qquad k=1,\ldots, K^{-i} 
\end{cases}
\end{equation*}
and $g_0(y_i) = \int g(y_i|\theta) dF_0(\theta)$. Intuitevely, $c_i = 0$ is the case where subject $i$ is assigned to a new atom in an occupied component, and $c_i=K^{-i}+1$ is the case where subject $i$ is assigned to an atom of a new component. 

Then the atoms are updated from 
$$
\theta_k | y,z,c \propto \left\{ \prod_{i: c_i = k} g(y_i|\theta_k) \right\} F_0(\theta_k)
$$
The variables $V_k$ are updated through data augmentation, by introducing $A_{ik} \sim Bern (V_k)$ and $B_{ik} \sim Bern(w(s,t,\psi,\zeta))$, independently for each $k$. Then, let $H_i = \min \{k: A_{ik} = B_{ik} = 1\}$, we have that 
$$
V_k \sim Be\left(a + \sum_{i: H_i \geq k} A_{ik}, b + \sum_{i: H_i \geq k} (1-A_{ik}) \right).
$$
Updating $\psi_k$ and $\zeta_k$ requires Metropolis-Hastings steps. 

Throughout the paper, we use mixtures of Gaussian distributions. Strategies to update the component specific parameters $(\mu_k,\sigma^2_k)$ are typical of mixture models, i.e. 
$$
\mu_k \sim \mathcal{N}\left(\bar{y}_{k}, \tau^2_{\mu} + n_k \sigma^2_k\right) \\
$$
where $\bar{y}_{k}$ is the average of the observations allocated to component $k$, $n_k$ is the number of observations allocated to component $k$, and $\tau^2_{\mu}$ is the prior variance (assuming conjugate prior distributions for the parameters). For $\sigma^2_k$, a Metropolis-Hastings step can be used. 

\section{Simulation study}
\label{sec:simu}

\subsection{Simulation study with a complex temporal structure}
\label{sub:simu_tempo}

We first illustrate the performance of the model following the simulation study set-up used by \cite{gutierrez2016time}, with the difference that we  include spatial dependence. 

We simulate $n=7,200$ observations for each time point $t=1,2,\ldots, 24$ from
$$
y(s,t) = 
\begin{cases}
\mathcal{N}(f(t),\sigma^2_1 C(s;\rho)) & 1 \leq t < 8 \\
0.3 \mathcal{N}(f(t),\sigma^2_1 C(s;\rho))) + 0.7 \mathcal{N}(f(t),\sigma^2_2 C(s;\rho))) & 8 \leq t  < 16 \\
0.5 \mathcal{N}(f(t),\sigma^2_3 C(s;\rho))) + 0.5 \mathcal{N}(0.1t + f(t),\sigma^2_3 C(s;\rho))) &  t \geq 16 \\
\end{cases}
$$
where $\sigma_1^2 = 0.04$, $\sigma^2_2=1$, $\sigma^2_3 = 0.09$ and
$$
f(t) = \cos(t) + 2 \times \sin(t) + \frac{t}{2} - \min(t,16).
$$
and $C(s;\rho))$ is a covariance matrix defined by a squared exponential kernel with lengthscale $\rho>0$. 
$\mathcal{N}(\mu, \Sigma)$ stands for a multivariate normal distribution, where $\Sigma$ is given an inverse-Wishart prior distribution with expected value $\mathbb{E}(\Sigma)=\frac{A}{\nu -3}$ and $\nu = 8$. We assume no covariates are available.

For this setting and the relative sample size, it is not possible to apply the version of the spatio-temporal stick-breaking with varying atoms of Section \ref{sub:varyingweights}, for reasons associated to both memory and computational time. Therefore it has been fitted by using only 50\% of the observations as training set. 

The results of the single-atom spatio-temporal stick-breaking (stSB), and its extension with varying atoms (stSB (VA)), are compared with two main alternatives available in the literature. The benchmark model is Bayesian spatio-temporal model (sp):
$$
y(s,t) = \theta(s,t) + \varepsilon(s,t)
$$
where $\theta(s,t)$ is a spatio-temporal process with a separable covariance function $\Sigma = C_{st} = C(y(s,t),y(s',t')) = C_s(s,s') \cdot C_t(t,t')$. $C_s(s,s')$ and $C_t(t,t')$ represent the spatial and temporal covariance functions, respectively. Then $\theta(s,t) \sim \mathcal{GP}(0,\tau^2 C_{st})$, where $\mathcal{GP}$ stands for a Gaussian process with a column vector of zero-elements as mean. Moreover, $\varepsilon(s,t)\sim \mathcal{N}(0,\sigma^2_{\varepsilon} \mathbf{I})$, where $\mathbf{I}$ is the identity matrix. The model can be fitted by using the \texttt{spTimer}  package \citep{bakar2015sptimer} available for the \texttt{R} software \citep{rcoreteam}. In particular, we choose to use an exponential correlation functions $C_s(s,s') = \exp(- \rho d(s,s')^2)$ and $C_t(t,t') = \exp(-\xi ( t - t')^2)$, where $d(s,s')$ is the Euclidean distance, and $\mathbb{V}\mbox{ar}(\theta(s,t)) =\tau^2$. Prior distributions are defined as $\rho, \xi \sim \mathcal{U}(3,30)$, $\sigma^2_{\varepsilon} \sim \mathcal{IG}(0.01,0.01)$, and $\tau^2 \sim \mathcal{IG}(0.01,0.01)$, where $\mathcal{IG}$ stands for an inverse gamma distribution.  

The second model is the one implemented by \cite{reich2007multivariate} (sSB), described in Section \ref{sec:spatialSB}, with prior distributions $a \sim Unif(0,10)$, $b\sim Unif(0,10)$, $\sigma^2_{\varepsilon} \sim \mathcal{IG}(0.01,0.01)$ and $\tau^2 \sim IG(0.01,0.01)$, the knots are uniform over the spatial domain; an exponential kernel is used for $w_k(s,\psi_k)$, with $h_{jk} \sim IG\left( 1.5, \frac{\nu^2}{2}\right)$ and $\nu \sim Unif(0, \nu_{max})$ where $\nu_{max}$ is the maximum distance between any pair of points in the spatial grid. Temporal dependence is introduced in the atoms of the stick-breaking process, by adding a temporal structure through a transition equation for $\theta(s,t)$, say $\theta(s,t) = \zeta \theta(s,t-1) + \eta(s,t)$, where $|\zeta| < 1$ and $\eta_{\mathcal{D},\mathcal{T}}=\{\eta(s,t): s\in \mathcal{D},t \in \mathcal{T}\}$ are independent realizations from a spatial stochastic process; see \cite{kottas2008modeling}.

We run the MCMC algorithm for 20,000 iterations, with the first 10,000 iterations discarded for burn-in. We set the maximum number of components to $M=100$. Data simulation is repeated 50 times. 

We use the expected mean squared prediction error (ESPE) to compare the prediction accuracy for each considered model; we keep 70\% of the simulated data as training dataset, and use the estimated model to predict the remaining 30\% of the observations. The ESPE is computed over the repetitions of the simulations and over all the time points and spatial locations:
$$
\mbox{ESPE} = \mathbb{E}\left[ \sum_{t \in \mathcal{T}} \sum_{s \in \mathcal{S}} [y(s,t) - \tilde{y}(s,t)]^2 \right],
$$
which is approximated by the corresponding sample average. Here, $y(s,t)$ is the observed value, while $\tilde{y}(s,t)$ is the predicted mean value. 

The mean squared error is 175.69 for sp, 196.04 for sSB, and 52.33 for stSB; for stSB, 97.3\% of accepted values of $\lambda$ are equal to zero and the posterior mean of values different from zero is 0.037. The stSB with varying atoms has ESPE equal to 142.15, with 100\% of accepted values of $\lambda$ equal to zero; however, it has to be noticed that this last model is fitted with a reduced training set with respect to the other models. By considering the prediction (not shown here), it is evident that the spatio-temporal model smooths too much over time, while the sSB is associated with the highest uncertainty in the prediction. The stSB is able to follow the data more closely. 

\subsection{Simulation study with a complex spatial structure}
\label{sub:simu_space}

We now consider the case where data are simulated from different spatio-temporal random fields with consecutive 15 time points and a variable number of locations which are simulated in the following way. Locations are generated according to a Thomas point process, i.e. a cluster point process, so that locations tend to get clustered: first a homogeneous Poisson point process with intensity $\omega > 0$ is simulated on a spatial rectangle (parent process), then for each simulated point another Poisson point process with intensity $\delta > 0$ is simulated (daughter processes). On these simulated daughter points, a random field is generated. The daughter point processes are simulated within a cluster disk (of radius 0.1), which has the centre in each parent point. This way of simulating spatial points has the aim of recreating a realistic situation where observations are locally closed; for instance, meteorological variables are often recorded by monitoring stations which are present in urban areas in a large number, while are less present or totally absent in rural areas. 
 
Some of the models include temporal dependence and some do not; this decision is meant to show the ability of the proposed approach to adapt to less complicated situations. We considered the following models:
\begin{enumerate}
	\item a stationary isotropic covariance model with Gaussian covariance 
	$$
	C(s,s') = \tau^2 e^{-\frac{||s-s'||^2}{h}}
	$$	
	with variance $\tau^2=1$, scale parameter $h=0.4$, and time-independent observations;
	\item a stationary isotropic exponential model
	$$
	C(s,s') = \tau^2 e^{-\frac{|s-s'|}{h}}
	$$	
	with variance $\tau^2=4$ and scale $h=10$, with a nugget component and a mean described by a trend of mean 0.5, and time-independent observations;
	\item a stationary isotropic covariance model belonging to the stable family, with covariance function given by 
	$$
	C(s,s') = \tau^2 e^{-\frac{|s-s'|^{\alpha}}{h}}
	$$
	with $\alpha=1.9$ so that the Gaussian sample paths are not differentiable, scale parameter $h=0.4$, and time-independent observations;
	\item a nugget effect model with variance $\tau^2=1$, zonal anisotropy, and time-independent observations;
	\item a Stein non-separable space-time model \citep{stein2005space}, where the covariance function is given by
$$
	C(h,t) = W_{\nu}(z) - \frac{(< h,c > t)}{((\nu-1)(2\nu+d)} \cdot W_{\nu-1}(z)),
	$$ 
where $h = ||s-s'||$, $W_{\nu}(\cdot)$ is a Whittle-Matern covariance model with smoothness parameter $\nu$, $z = ||(h,t)||$ is the norm of the vector $(h,t)$, and $d$ is the dimension of the space on which the random field is defined; here $\nu=1.5$; $t$ is intended as the difference between the observed time and the starting point and $<\cdot>$ is the inner product with respect to a reference vector of coordinates $c = (0.9,0.1)$;
	\item a stationary space-isotropic covariance model with covariance function given by
	$$
	C(s,s',t,t') = (\psi(t,t')+1)^{d/2} \phi\left( \frac{|s-s'|}{\sqrt{\psi(t,t')+1}}\right)
	$$
	where $d$ depends on the spatial dimension of the field and here is set equal to 2; $\phi(\cdot)$ is a Gaussian covariance model; $\psi(t,t') = |t-t'|^{\alpha}$ with $\alpha=1$.
\end{enumerate}

Simulations are obtained via the \texttt{R} package \texttt{RandomFields} \citep{schlather2015analysis}. For every model, we repeat simulations 50 times. We compare two exponential models (Model 2) with nugget and trend, using $\omega=10$ (Model 2a) and $\omega = 30$ (Model 2b) respectively in the Thomas point process.

Results are shown in Table \ref{tab:simu}. The stSB outperforms, in terms of ESPE, sSB for simulations from a Gaussian model (Model 1), an anysotropy model (Model 4), a Stein model (Model 5), and a non-separable space-time model (Model 6). Results are very similar for the sp, sSB and stSB for simulations from exponential models (Models 2a and 2b). 
Overall, stSB seems to be a good prior process in more complex situations, where space and time are not separable or where the covariance function is complicated, as for the Stein model. The approach still performs well in comparison with a simpler kriging model in more standard situations. In this sense, stSB seems to represent a robust prior process. The stSB with varying weights and atoms tend to perform slightly worse than the other processes with simulations from the exponential models (Model 2a and 2b), the stable model, the Stein model and the non-separable model. This may be due to the additional complexity of the model, which requires an additional number of MCMC iterations to provide good estimates and an additional number of observations to reach the same level of accuracy. 

\begin{table}[]
\centering
\begin{tabular}{l|cccc}
                                    & \textbf{sp} & \textbf{sSB} & \textbf{stSB} & \textbf{stSB (VA)}\\ \hline
\textit{1. Gaussian}                   & 0.93611              & 0.96797 & 0.92212      & \textbf{0.89518}          \\
\textit{2a. Trend, $\omega=10$}               & \textbf{1.77601}               & 1.78522 & 1.86372 & 2.86420                \\
\textit{2b. Trend, $\omega=30$} & \textbf{1.66367}               & 1.68756 & 1.71742        &     1.79734  \\
\textit{3. Stable}                     & \textbf{0.11978}              & 0.31019 & 0.90373   & 0.98556      \\
\textit{4. Zonal Anisotropy}                & 1.00187              & 0.99860 & \textbf{0.90475}         &    0.98629   \\
\textit{5. Stein}     &  0.63084 &   0.68333           &  \textbf{0.57722}      & 0.77151               \\
\textit{6. Nonseparable}              & 0.60420             & 0.64485 & \textbf{0.53416} & 0.93089
\end{tabular}
\caption{Comparison of the estimated ESPEs for the four considered models and for simulations from models described in the list with the relative number. The bold number represents the preferred models.}
\label{tab:simu}
\end{table}

\begin{table}[]
\begin{tabular}{l|cccc}
                                      & \textbf{stSB}                          & \textbf{stSB (VA)}                     & \textbf{stSB}                  & \textbf{stSB (VA)}             \\
                                      & $\Pr(\lambda=0 | \mathbf{y})$ & $\Pr(\lambda=0 | \mathbf{y})$ & $\mathbb{E}[\lambda | \mathbf{y}]$ & $\mathbb{E}[\lambda | \mathbf{y}]$ \\ \hline
\textit{1. Gaussian} & 0.97476                     & 0.00000                         & 0.04098            & 0.86898             \\
\textit{2a. Trend, $\omega=10$}                                                &     0.70008                         &                              0.15749 &    0.34658                 &                      0.94788 \\
\textit{2a. Trend, $\omega=30$}                                      &                              0.69074 &   0.50000                            &                      0.33236 &    0.99458                   \\
\textit{3. Stable}                                      &   0.87393                            &   0.10000                            &                      0.13651 &    0.91893                   \\
\textit{4. Zonal Anisotropy}                                                     &       0.02957                        &   0.08295                            &       0.94224                &         0.89265              \\
\textit{5. Stein}                                     &     0.98938 &    0.10000                           &                      0.01257 &     0.86846                  
\\
\textit{6. Nonseparable}                                  &                              0.98252 &         0.00000                      &                      0.03513 &      0.97822                
\end{tabular}
\caption{Comparison of posterior measures on the separability parameter $\lambda$ in the spatio-temporal stick-breaking models - posterior probability that $\lambda=0$ and posterior expectation of $\lambda$ when different from zero.}
\label{tab:lambdatab}
\end{table}

Table \ref{tab:lambdatab} presents results on the estimation of the separability parameter $\lambda$ in Equation \eqref{eq:gneiting}. When data comes from a separable model, in particular, when observations are time-independent, there is a strong support that weights generated with a separable kernel are preferred. It is, however, interesting to notice, that the support for a non-separable kernel is never strong, except for the anisotropic case: data are evidently not strongly informative about the parameter defining separability within the kernel. This is not surprising, since the covariance among observations obtained by using either a separable kernel or a non-separable kernel is non-separable in both cases. However, it represents an additional level of flexibility for the model, allowing for the covariance function to be non-stationary. 
 
\begin{figure}[h]
	\includegraphics[width=10cm,height=8.5cm]{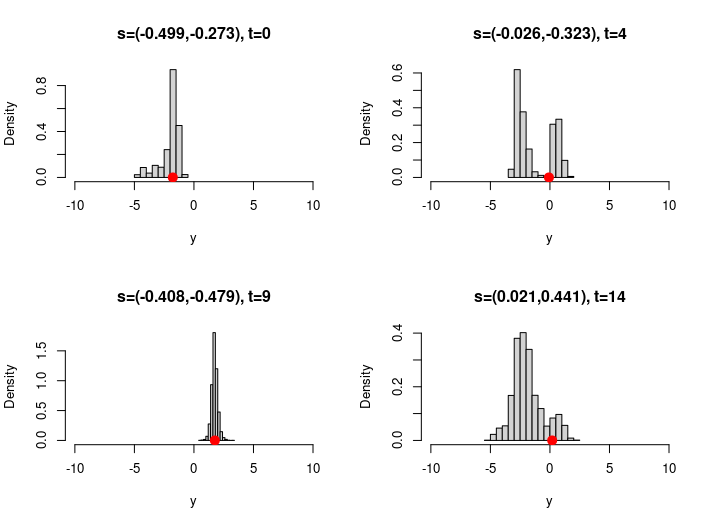}
	\caption{Mixture densities obtained by choosing the modal means $\mu$, weights $\pi$ and standard deviation $\sigma$ parameters for four observations, characterized by different spatial locations and time points and for data generated from a Stein model. The red dots on the $x$-axis represent the observed value.}
	\label{fig:stein_post}
\end{figure}

From the simulation study, it seems that stSB is able to adapt to spatial and temporal changes, by modifying the mixing weights according to the information carried by the data. Figure \ref{fig:stein_post} shows the mixture densities obtained via posterior estimate of the parameters relative to four locations and time points from data generated from a Stein model: the process seems to be adaptable to different situations and produce very variable posterior densities. 

\section{Real-data application}
\label{sec:real}

We now implement the methodology proposed in this work on a real-data application, the Supplementary Material includes an alternative example on temperature data. 

\subsection{Rainfall data}

Rainfall in Australia are known to be highly variable \citep{nicholls1997australian}, with low rainfall over most of the country and intense seasonal falls in the tropical regions in the North. The rainfall pattern is concentric from the coast, in particular the North and East coasts can be characterized by large rainfall volumes, while the centre is more arid. Around the centre there is another humid area, particularly in the East \citep{dey2019review}. The temporal and spatial variability of Australia's rainfall can lead to droughts and floods that have large impact on various sectors, e.g. agriculture, hydrology, resource management and allocation, human health and housing. 

\begin{figure}[h]
\begin{multicols}{3}
    \includegraphics[width=\linewidth]{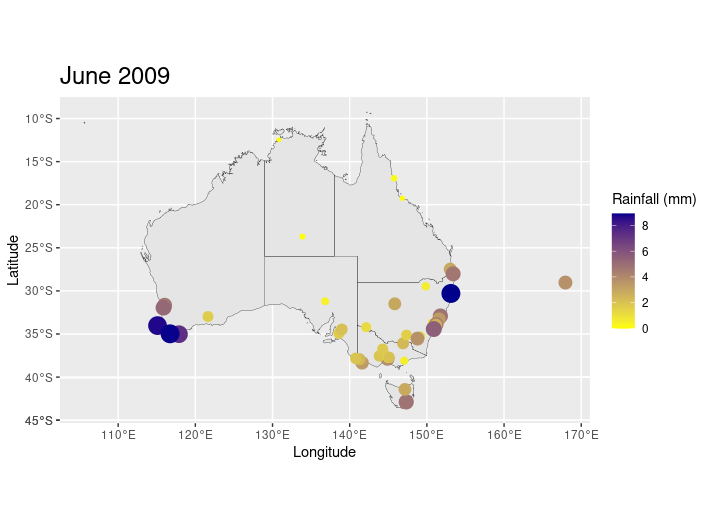}\par 
    \includegraphics[width=\linewidth]{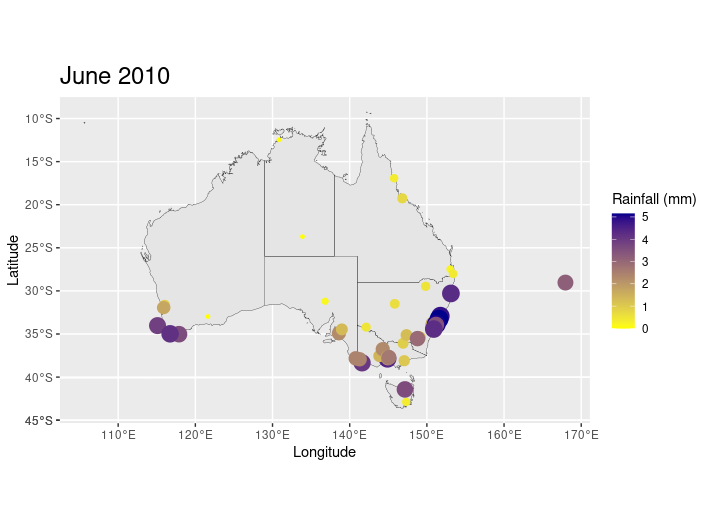}\par 
    \includegraphics[width=\linewidth]{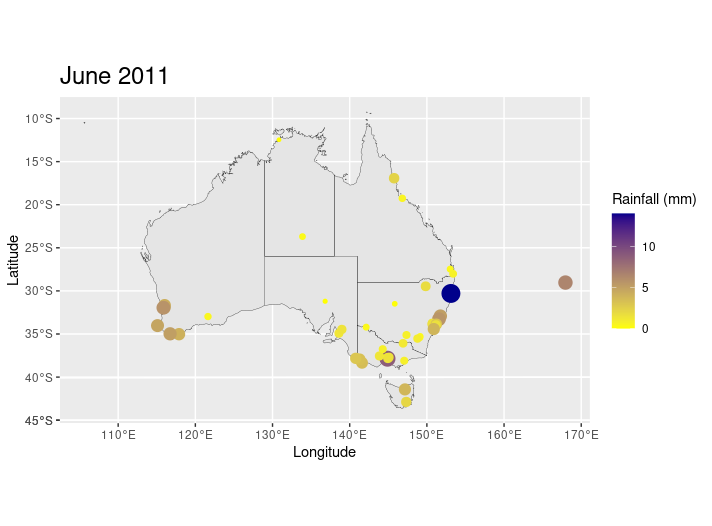}\par 
    \end{multicols}
\begin{multicols}{3}
    \includegraphics[width=\linewidth]{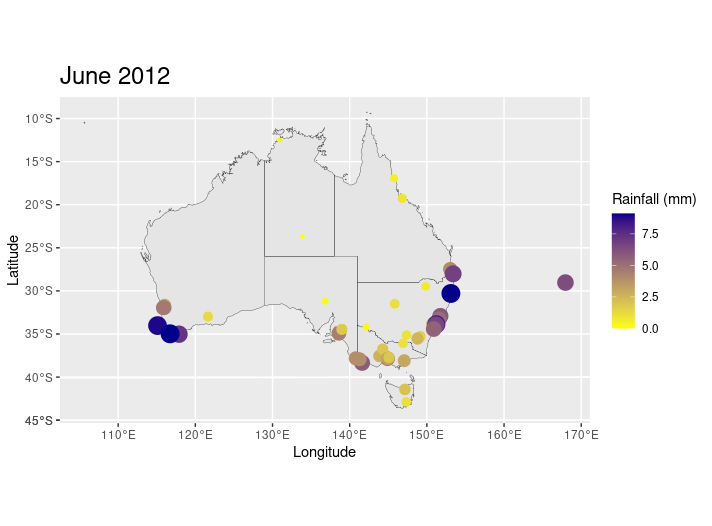}\par 
    \includegraphics[width=\linewidth]{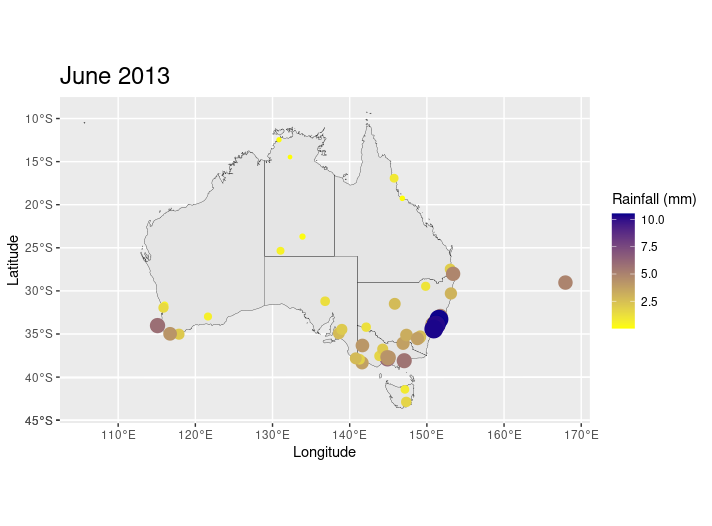}\par 
    \includegraphics[width=\linewidth]{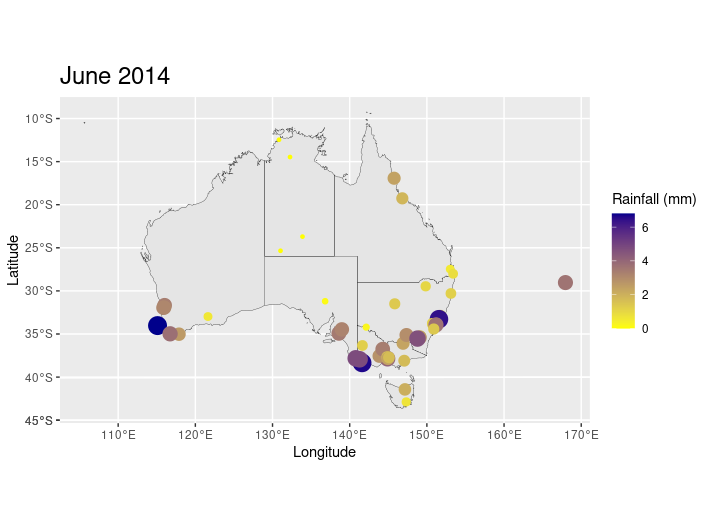}\par 
    \end{multicols}
\begin{multicols}{3}
    \includegraphics[width=\linewidth]{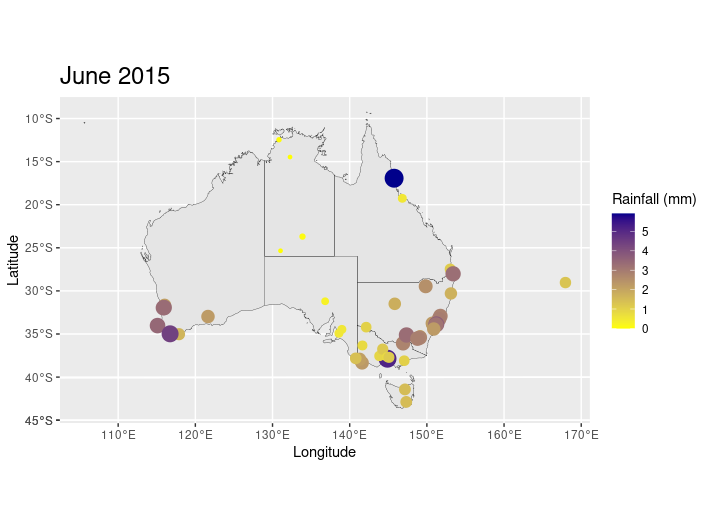}\par 
    \includegraphics[width=\linewidth]{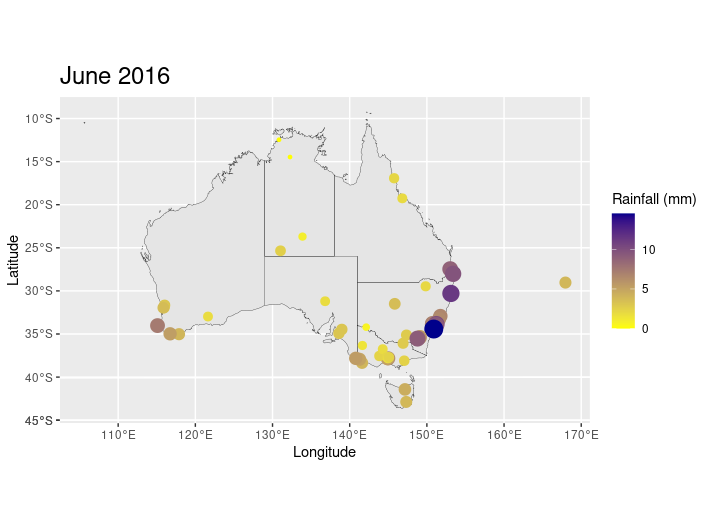}\par 
    \includegraphics[width=\linewidth]{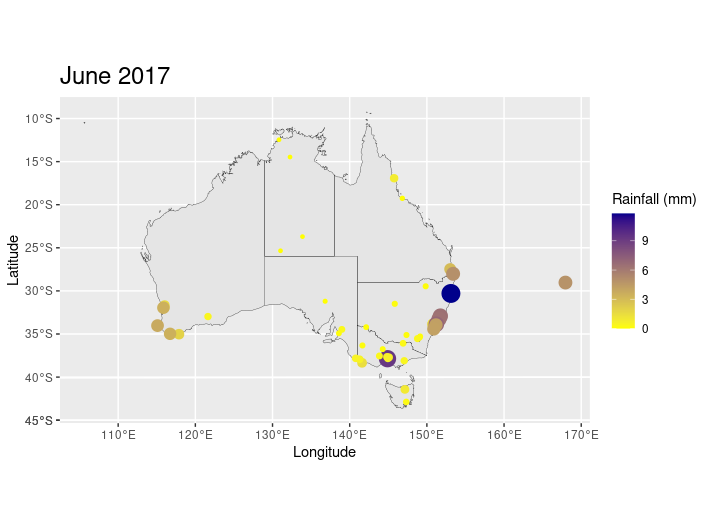}\par 
    \end{multicols}
\caption{Monthly rainfalls for the month of June from 2009 to 2017 in Australia. A larger and darker dot represents a larger volume of rainfall.}
\label{fig:rainfall_june}
\end{figure}

The dataset contains 10 years of daily weather observations in 49 locations across all Australia, from 31 October 2007 to 24 June 2017. Observations were drawn from weather stations, and include the amount of rainfall recorded for the day in mm. The dataset is publicly available on the website of the Bureau of Meteorology of the Australian Government \url{http://www.bom.gov.au/climate/data.} For illustration, Figure \ref{fig:rainfall_june} shows the monthly rainfall volumes for each location from 2009 to 2017. Years 2007 and 2008 were left our of the illustration because data were highly missing, however they were included in the training dataset. 

Four models have been compared: sp, sSB, single-atom stSB and stSB with varying atoms and weights. The first three models were trained on the daily observations from 2007 to 2016 to predict data for 2017; stSB (VA) could not be fit on the full dataset, for computational reasons. Therefore, it was fitted on the weakly data, instead of daily ones, obtained by avaraging the daily rainfall data for each week at each location. Models were compared by looking at the ESPE: 109.74 for sp, 112.09 for sSB, 71.05 for single-atom stSB, and 104.67 for stSB (VA). Again, stSB seems to produce more realible predictions with respect to the other models. A fair comparison with the most complicated model, stSB (VA), is not possible, because the model was trained on a relatively different dataset.

\begin{figure}
\centering
\begin{subfigure}{.33\textwidth}
  \centering
  \includegraphics[width=5cm,height=4.2cm]{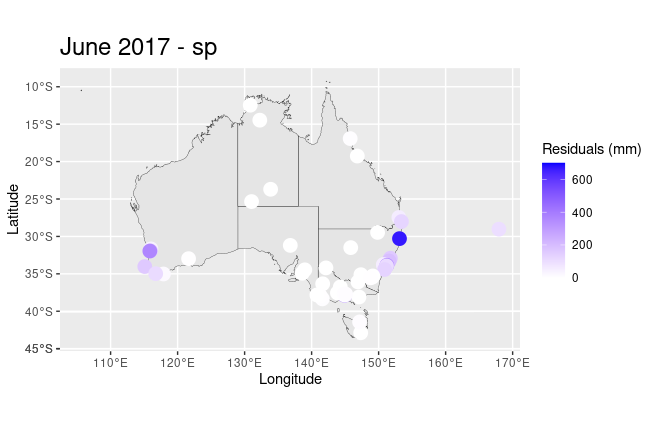}
\end{subfigure}%
\begin{subfigure}{.33\textwidth}
  \centering
  \includegraphics[width=5cm,height=4.2cm]{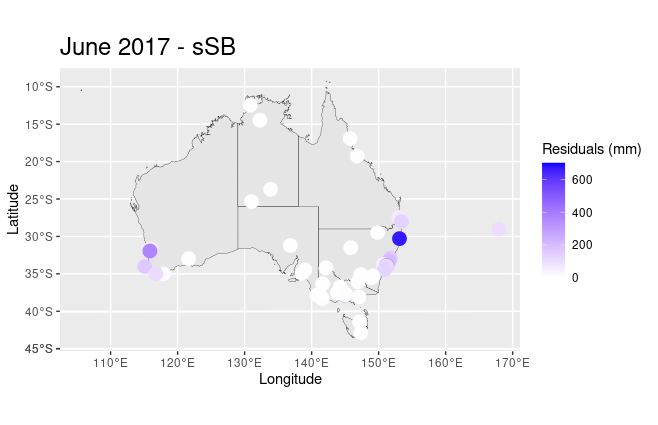}
\end{subfigure}
\begin{subfigure}{.33\textwidth}
  \centering
  \includegraphics[width=5cm,height=4.2cm]{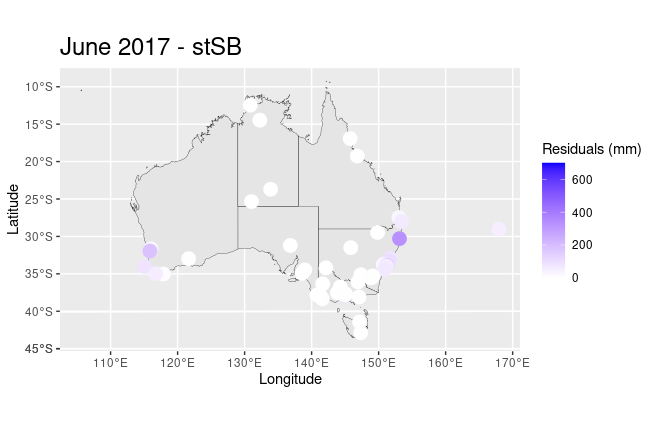}
\end{subfigure}
\caption{Residuals plots for the month of June 2017 obtained with the three methods, sp (left), sSB (centre), stSB (right). The size of dots is fixed.}
\label{fig:res_rainfall}
\end{figure}

Figure \ref{fig:res_rainfall} shows the average squared residuals for predictions relative to June 2017 (values based on the posterior predictive mean) for each city. The first two models, sp and sSB, tend to oversmooth over the area, with high volatility, in particular in the area relative to Perth (West) and the area on the coast of New South Wales (East), especially around Byron Bay. On the other hand, stSB is able to capture the peaks of rainfall and the squared residuals show less spatial structure. 

For this application, the percentage of posterior samples for the separability parameter $\lambda$ which are equal to zero is 100\%: in this case, the hypothesis of separability between time and space for the mixing weights is not rejected. 

\section{Conclusion}
\label{sec:conclu}

Modelling spatio-temporal phenomena is an important and challenging problem. Sometimes there is not enough information to define a parametric model and nonparametric methods may be preferred. This work presents a semi-parametric method to model spatio-temporal data, which can be used to obtain higher predictive accuracy without the need of introducing a large number of covariates. 

The spatio-temporal stick-breaking prior proposed in this paper is a generalisation of \cite{reich2007multivariate} and \cite{dunson2008kernel}. This model avoids oversmoothing, but differently from the spatial stick-breaking of \cite{reich2007multivariate}, can model non-stationary covariance structures (with respect to either time or space, or both of them) and with a covariance among observations that can be symmetric or asymmetric with respect to space and time. Moreover, this model avoids the need for space-time varying covariates as in \cite{hossain2013space} to model the weight components.

Two models have been considered: a single-atom model where only the mixing weights are dependent on time and space, and a model where both atoms and weights are allowed to vary over time and space. The two models have shown to often lead to similar results, however the latter is definitely more computationally intensive, and this can lead to the need to train the model on a reduced number of observations. 

The model has been validated by using simulated data and real data, and compared to other spatio-temporal models. Two dataset has been analyzed: a dataset including temperature recordings between 1985 and 2004 in California, USA (available in the Supplementary Material), and a dataset of daily rainfall measurements in Australia. With respect to the spatio-temporal kriging model, the spatial stick-breaking and the spatio-temporal stick-breaking avoid oversmoothing, and with respect to the spatial stick-breaking, the spatio-temporal stick-breaking allows to better differentiate among areas over time: for example, predicting larger volumes of rain in areas that have become more humid in recent years, and which used to be drier in the past. 

An interesting current line of research focuses on situations where samples share some instead of all mixture components (with or without differing weights) and the components themselves may show some misalignment. For example,  \cite{lopes2003bayesian} and \cite{muller2004method} allows some clusters to be shared among samples, \cite{teh2004sharing} and \cite{cron2013hierarchical} allows all clusters to vary across samples, \cite{rodriguez2008nested} allows a structures where clusters in each sample must be either identical or completely different, and \cite{soriano2019mixture} propose a hierarchical mixture modelling framework allowing for both varying weights across samples and distribution misalignemnt. Such proposal often relies on a stick-breaking construction, however they also tend to assume independence within samples. The possibility to introduce spatial or temporal dependence on the mixing weights in a hierarchical way is left for further research. 

The model is open to several generalisations. First, we have intentionally considered cases where no covariates were available, but it is easy to incorporate covariate information. Similarly, it is possible to introduce seasonality components in the regression part of the model. Moreover, the model is easily generalisable to a multivariate setting.

\section*{Appendix A: The correlation structure among observations}

\begin{figure}[h]
	\includegraphics[width=10cm,height=7cm]{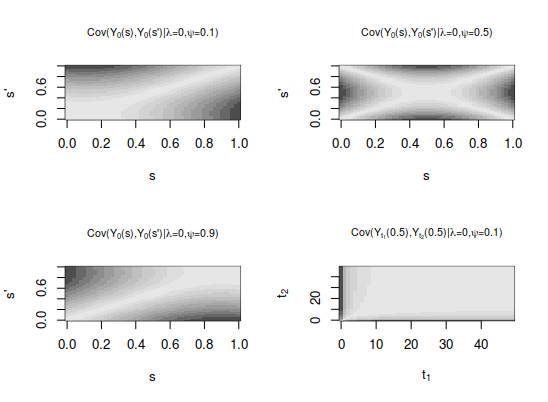}
	\caption{Correlation among observations at different spatial locations  with $\lambda=0$, and different choices of knots. The bottom right panel represents the correlation among observations at different time points for $s=s'=0.5$, $\lambda=0$ and $\psi=0.1$.}
	\label{fig:corr_lambda0}
\end{figure}

\begin{figure}
	\includegraphics[width=12cm,height=8cm]{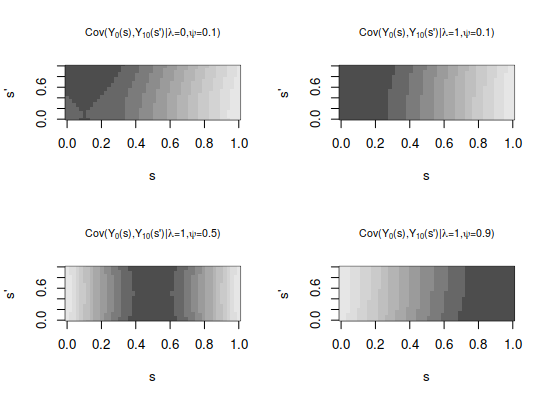}
	\caption{Representation of the correlation among observations at different spatial locations and at time $t=0$ and $t'=10$. The top left panel represent the case of separability between time and space ($\lambda=0$), while the other panels plots the cases of spatio-temporal interaction ($\lambda=1$), with $\psi=(0.1,0.5,0.9)$ respectively.}
	\label{fig:corr_lambda1}
\end{figure}

Figure \ref{fig:corr_lambda0} and \ref{fig:corr_lambda1} show examples of correlations among observations varying space or time, with $\gamma=1$ and for the extreme cases $\lambda=0$ and $\lambda=1$ and different choices of knots. Variables observed at the same time points show a correlation which is stronger as $s \rightarrow s'$; Figure \ref{fig:corr_lambda0} also shows that the correlation among observations at the same location is low when $t$ is far from $t'$, while it increases as $t,t' \rightarrow \infty$. Figure \ref{fig:corr_lambda1} analyzes the interaction between time and space, showing the correlation among observations recorded at different locations and different times. As the observations are recorded at different times, the correlation remains low where the knots are located. 

\newpage

\section*{Appendix B: Correlation among mixing weights}

\begin{figure}[h]
\includegraphics[width=12cm,height=12cm]{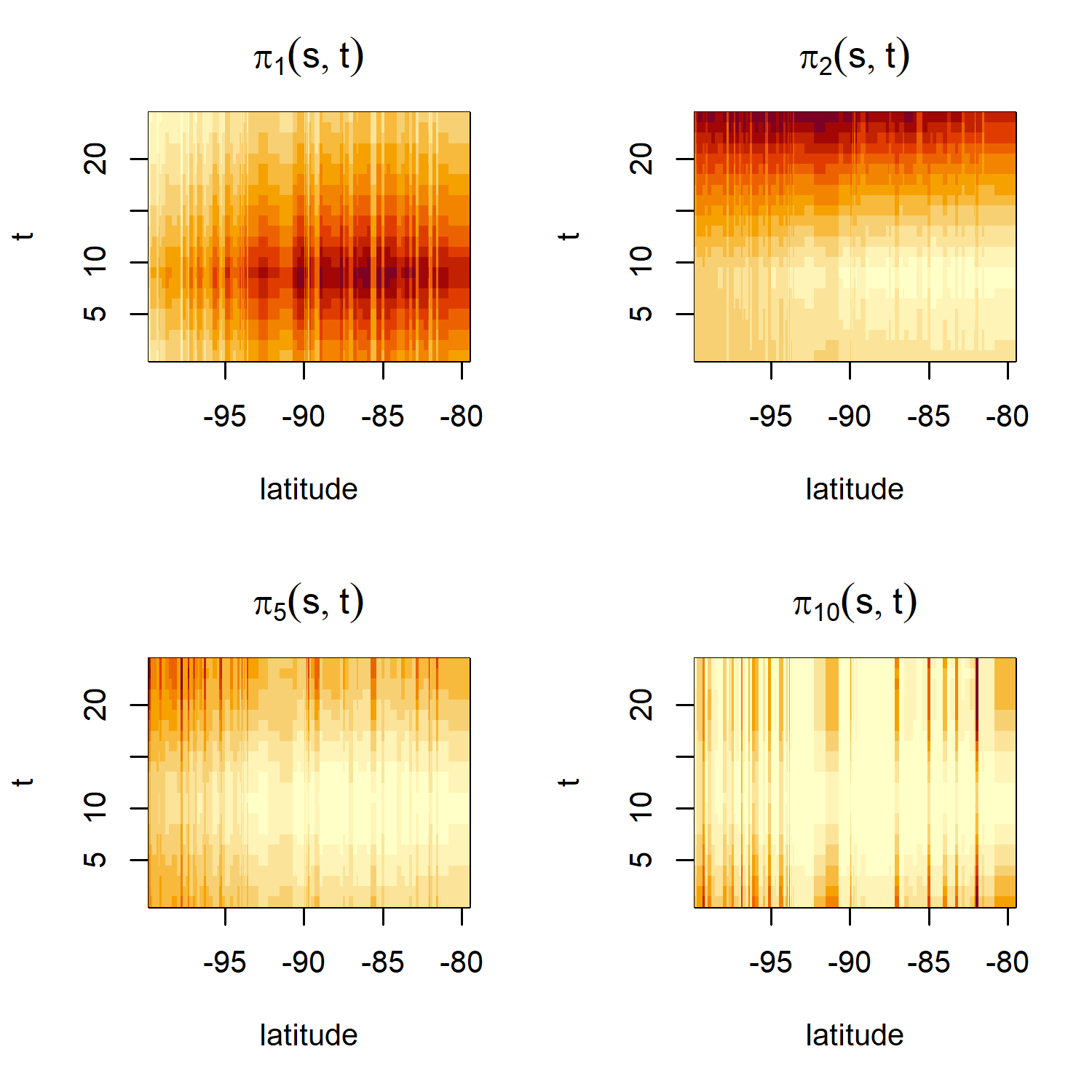}
	\caption{Contour plots of $\pi_1(s,t)$,$\pi_2(s,t)$,$\pi_5(s,t)$, $\pi_{10}(s,t)$ as a function of latitude and time.}
	\label{fig:pi_contour}
\end{figure}

Figure \ref{fig:pi_contour} shows the heatmaps of the variable $\pi(s,t)$ for different components and as a function of latitude and time; it shows how the mixing weights can strongly vary depending on time and space, in an asymmetric way. 

\pagebreak

\section*{Appendix C: Simulations from stSB and stSB (VA)}

\begin{figure}[h]
\begin{multicols}{2}
    \includegraphics[width=4.5cm,height=3.5cm]{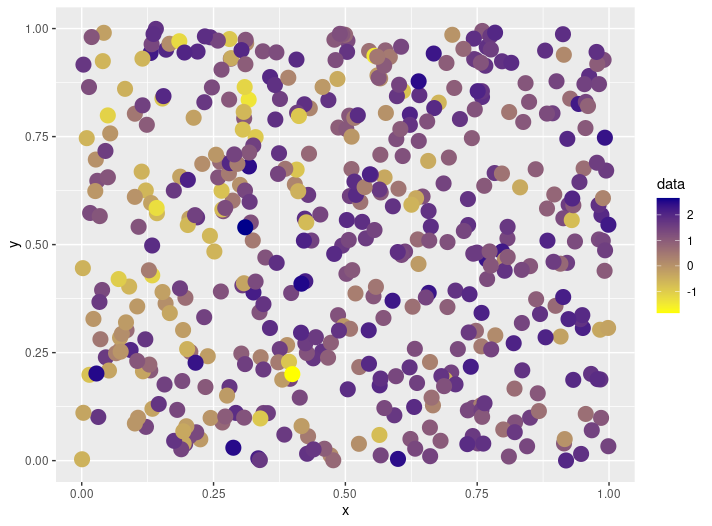}\par 
    \includegraphics[width=4.5cm,height=3.5cm]{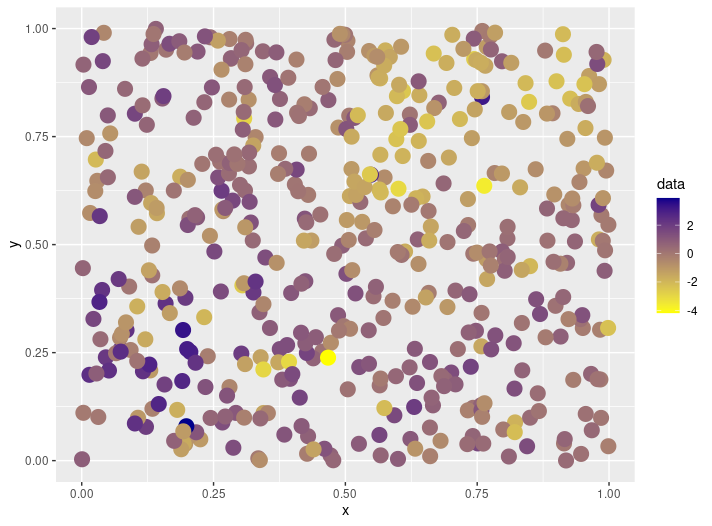}\par 
    \end{multicols}
\begin{multicols}{2}
    \includegraphics[width=4.5cm,height=3.5cm]{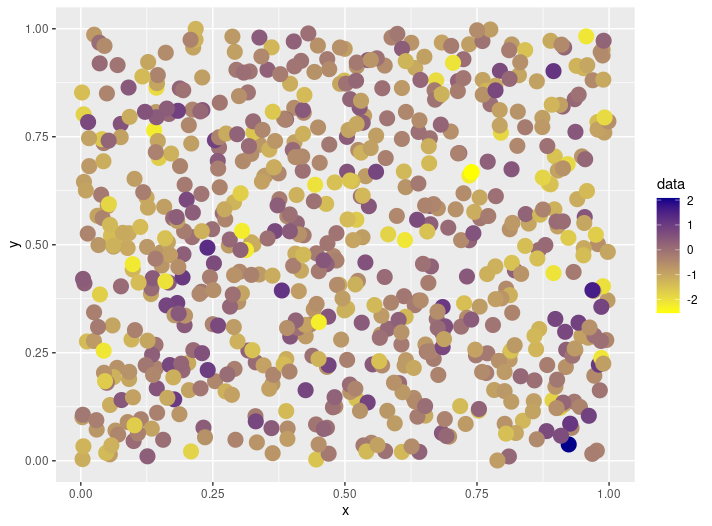}\par
    \includegraphics[width=4.5cm,height=3.5cm]{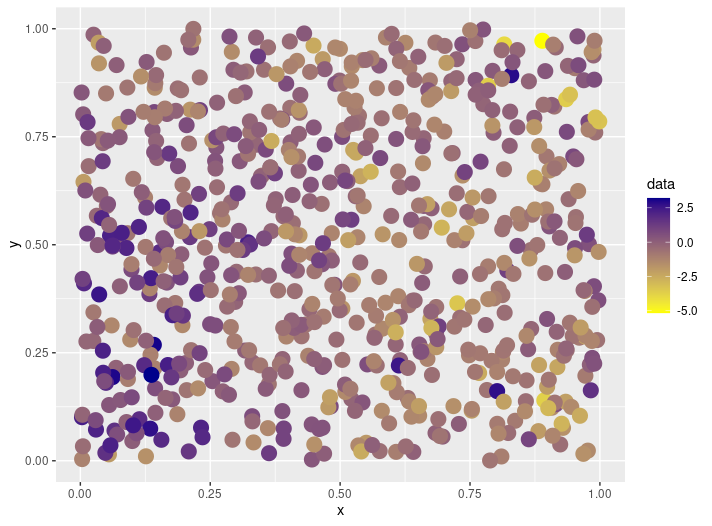}\par
\end{multicols}
\begin{multicols}{2}
    \includegraphics[width=4.5cm,height=3.5cm]{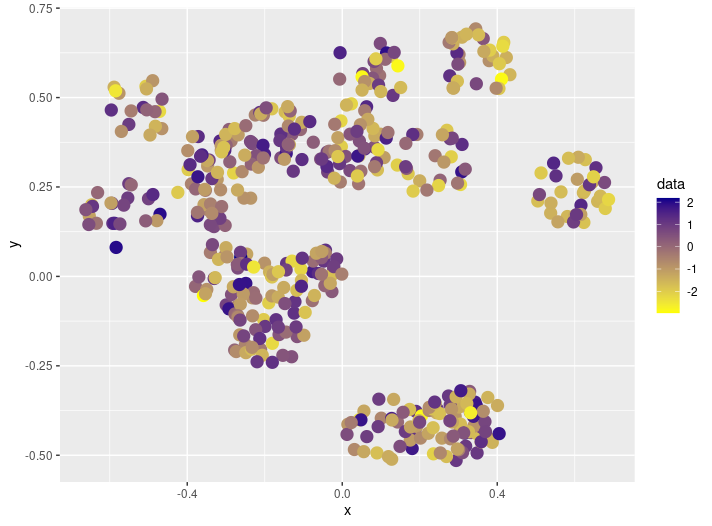}\par
    \includegraphics[width=4.5cm,height=3.5cm]{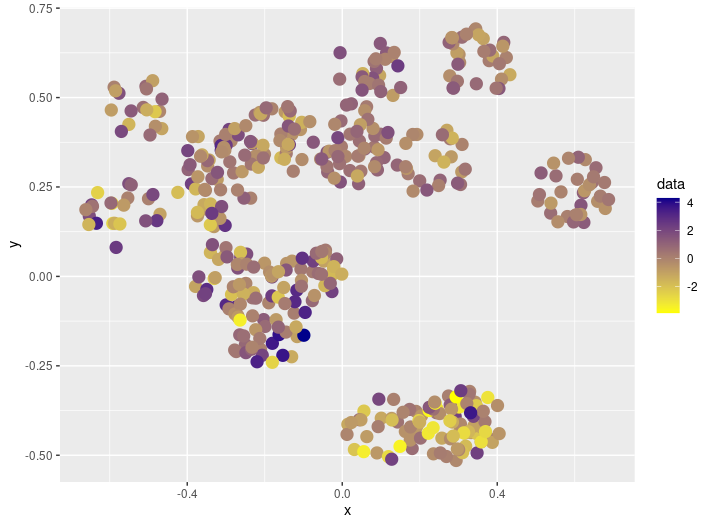}\par
\end{multicols}
\caption{Three examples processes simulated from a spatio-temporal stick-breaking with only weights depending on space and time (left) and with both weights and atoms depending on space and time (right). Top figure: 500 spatial locations simulated uniformly, $t=1$; middle figure: 700 spatial locations simulated uniformaly, $t=10$; bottom figure: 482 locations simulated according to a Thomas process, $t=1$.}
\label{fig:simulations}
\end{figure}

Figure \ref{fig:simulations} shows simulations from the spatio-temporal stick-breaking process with only weights dependent on time and space (left) and the spatio-temporal stick-breaking process with both weights and atoms dependent on time and space (right). The two processes seems to produce similar levels of dependence among observations, with the process introducing dependence also among atoms producing smoother variation over space, while the process allowing only the weights to be dependent on space and time seems to have more heterogeneity. 

\newpage

\begin{figure}[h]
\begin{multicols}{4}
	\textbf{$t=1$}\par
    \includegraphics[width=0.7\linewidth]{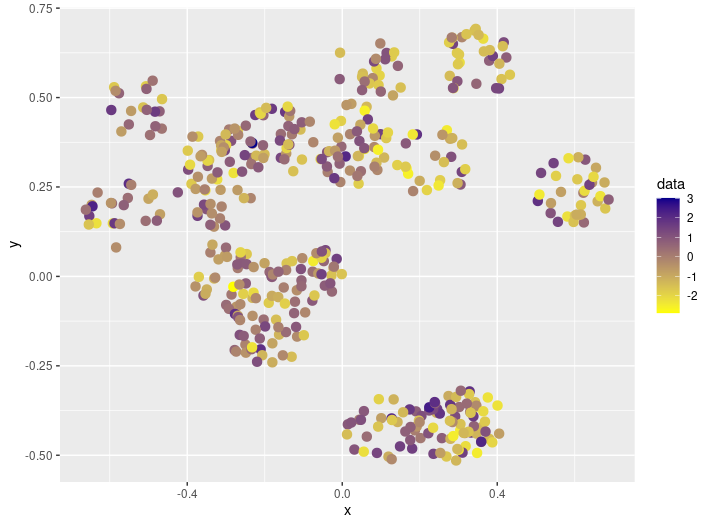}\par 
    \textbf{$t=2$}\par
    \includegraphics[width=0.7\linewidth]{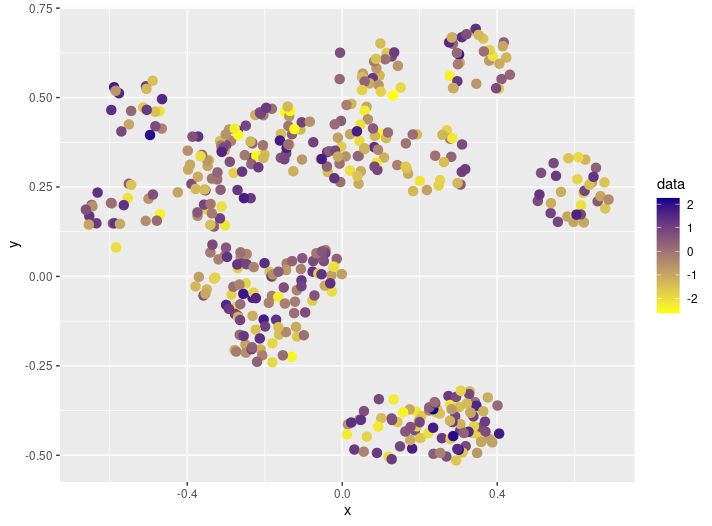}\par 
    \textbf{$t=1$}\par
    \includegraphics[width=0.7\linewidth]{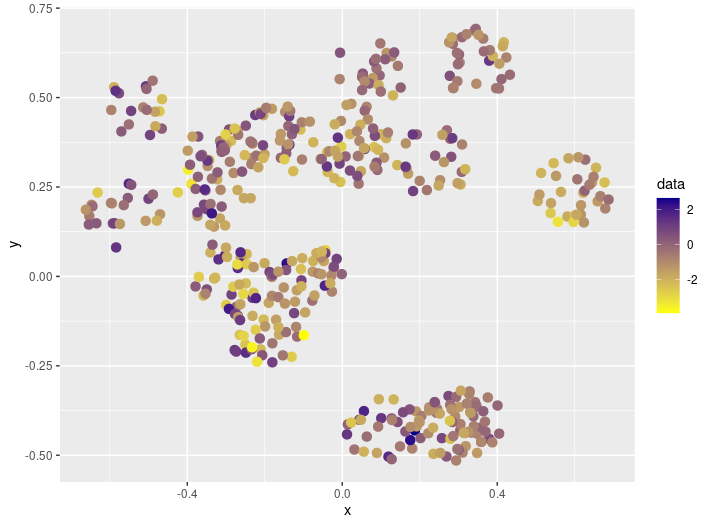}\par 
    \textbf{$t=2$}\par
    \includegraphics[width=0.7\linewidth]{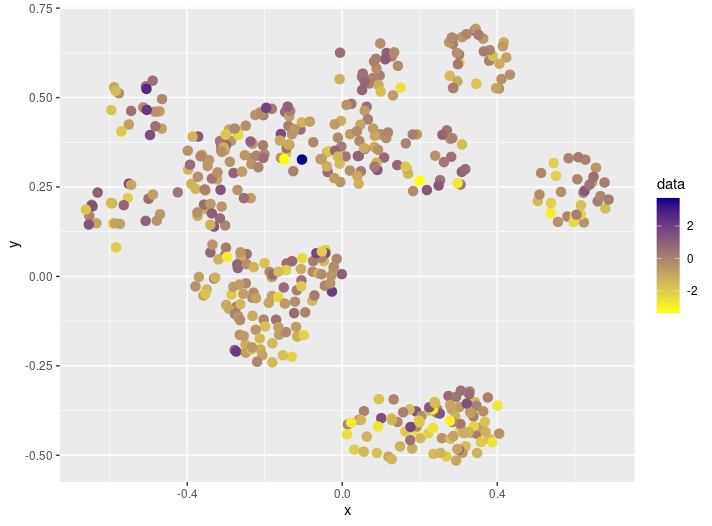}\par 
    \end{multicols}
\begin{multicols}{4}
    \textbf{$t=3$}\par
    \includegraphics[width=0.7\linewidth]{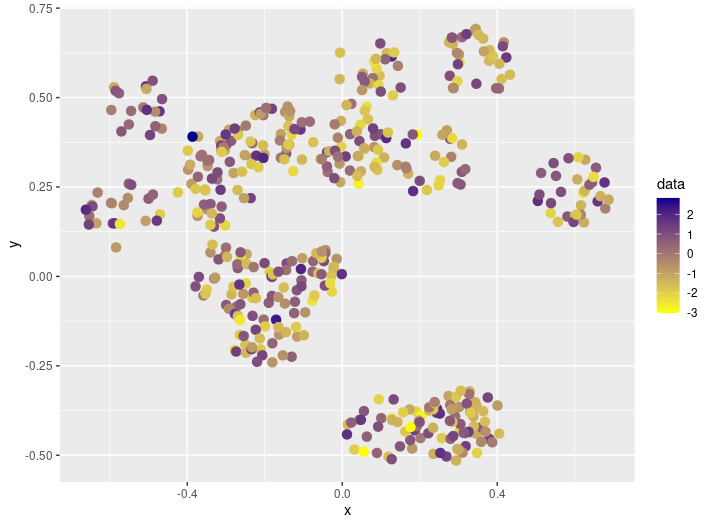}\par 
	\textbf{$t=4$}\par
    \includegraphics[width=0.7\linewidth]{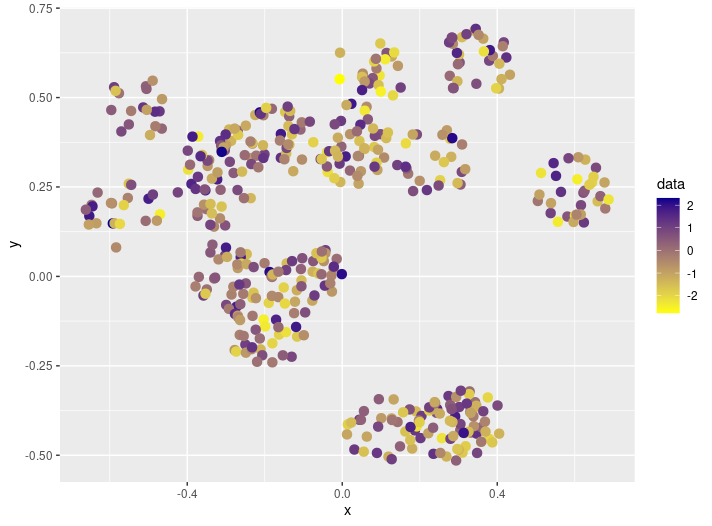}\par 
    \textbf{$t=3$}\par
    \includegraphics[width=0.7\linewidth]{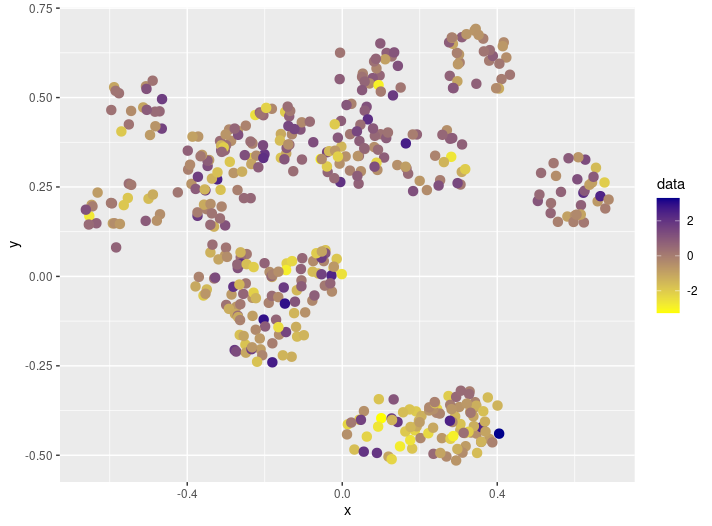}\par 
    \textbf{$t=4$}\par
    \includegraphics[width=0.7\linewidth]{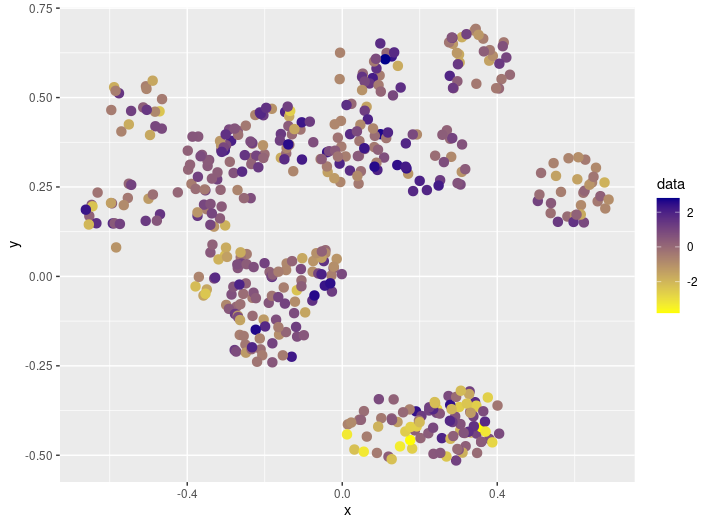}\par 
    \end{multicols}
\begin{multicols}{4}
    \textbf{$t=5$}\par
    \includegraphics[width=0.7\linewidth]{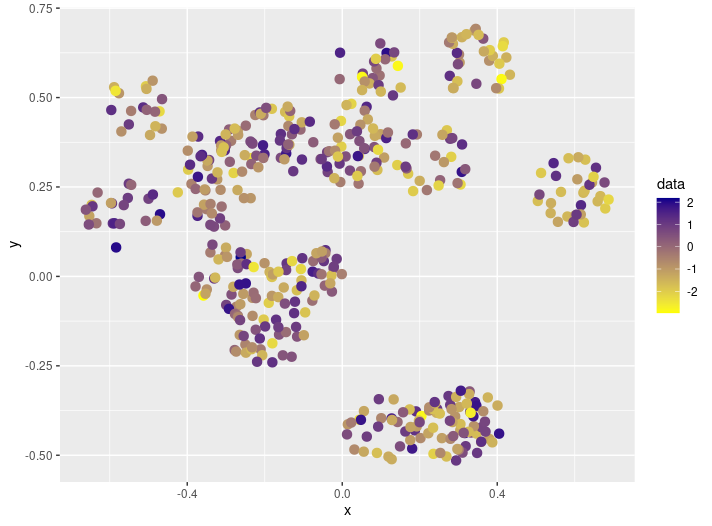}\par 
    \textbf{$t=6$}\par
    \includegraphics[width=0.7\linewidth]{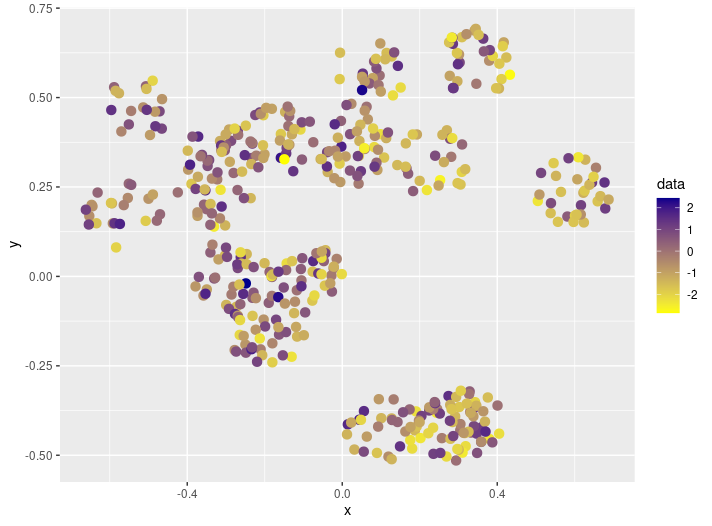}\par 
    \textbf{$t=5$}\par
    \includegraphics[width=0.7\linewidth]{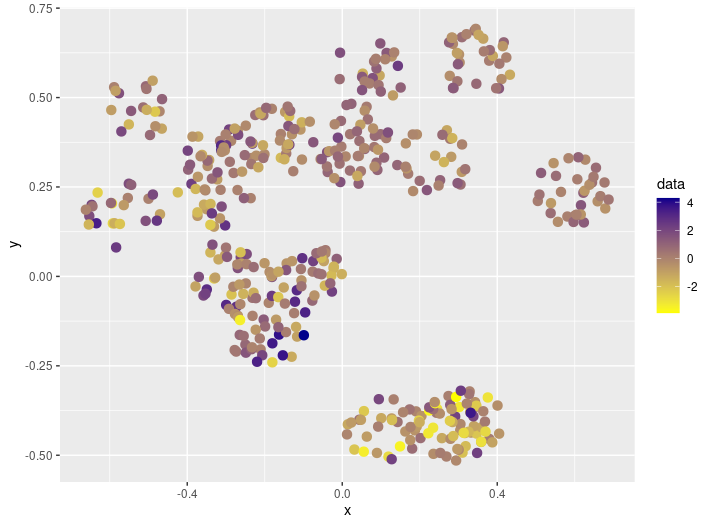}\par 
    \textbf{$t=6$}\par
    \includegraphics[width=0.7\linewidth]{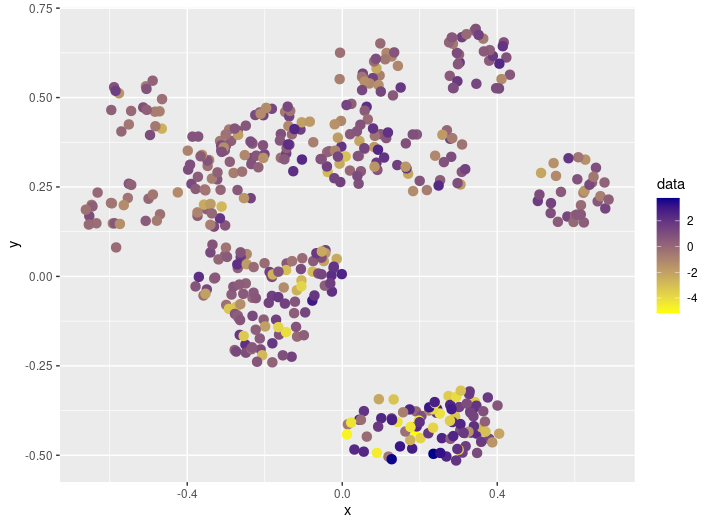}\par 
    \end{multicols}
\begin{multicols}{4}
    \textbf{$t=7$}\par
    \includegraphics[width=0.7\linewidth]{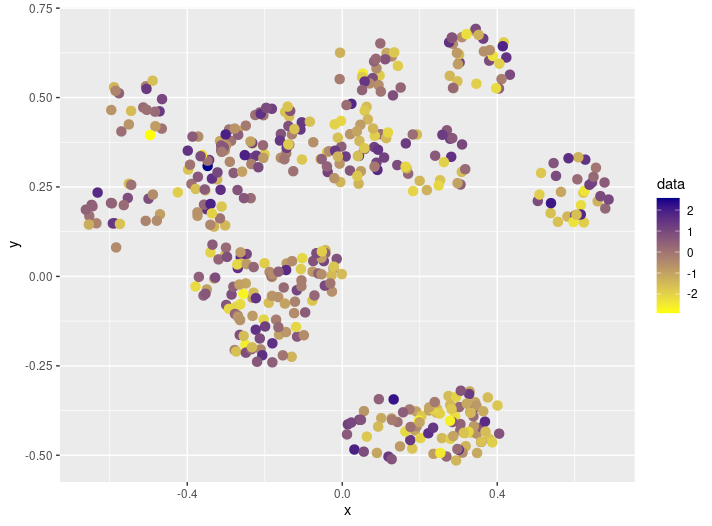}\par 
    \textbf{$t=8$}\par
    \includegraphics[width=0.7\linewidth]{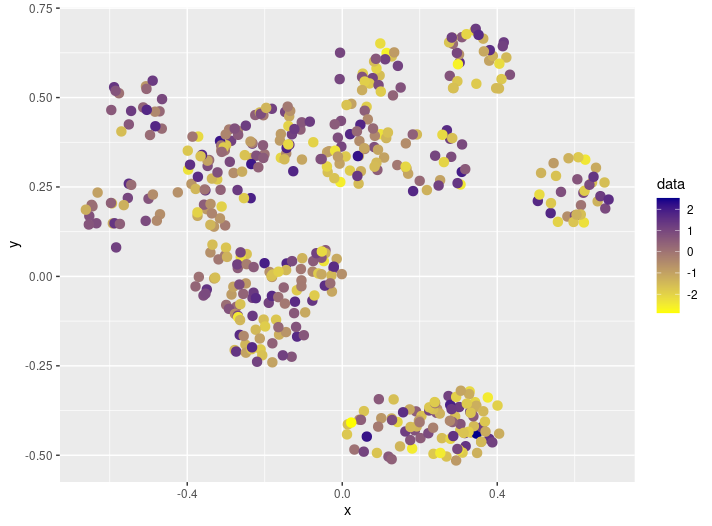}\par 
    \textbf{$t=7$}\par
    \includegraphics[width=0.7\linewidth]{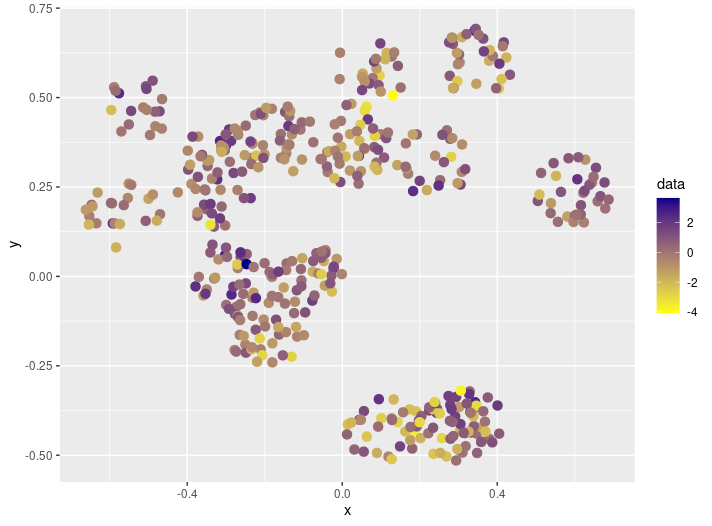}\par 
    \textbf{$t=8$}\par
    \includegraphics[width=0.7\linewidth]{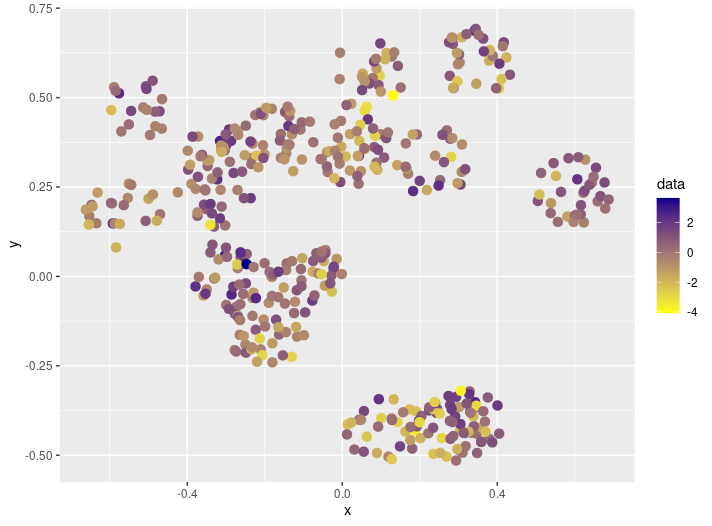}\par 
    \end{multicols}
\begin{multicols}{4}
    \textbf{$t=9$}\par
    \includegraphics[width=0.6\linewidth]{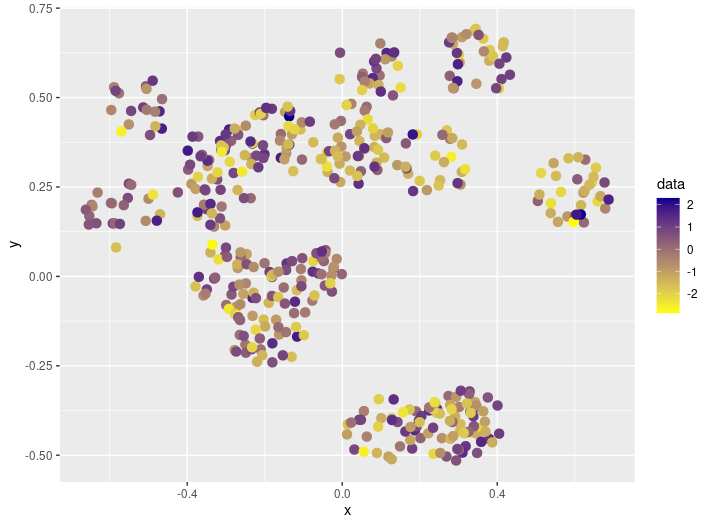}\par 
    \textbf{$t=10$}\par
    \includegraphics[width=0.6\linewidth]{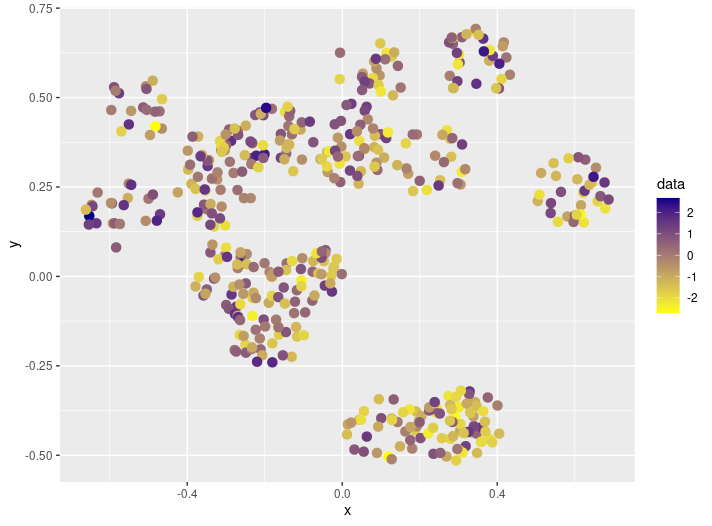}\par
    \textbf{$t=9$}\par
    \includegraphics[width=0.6\linewidth]{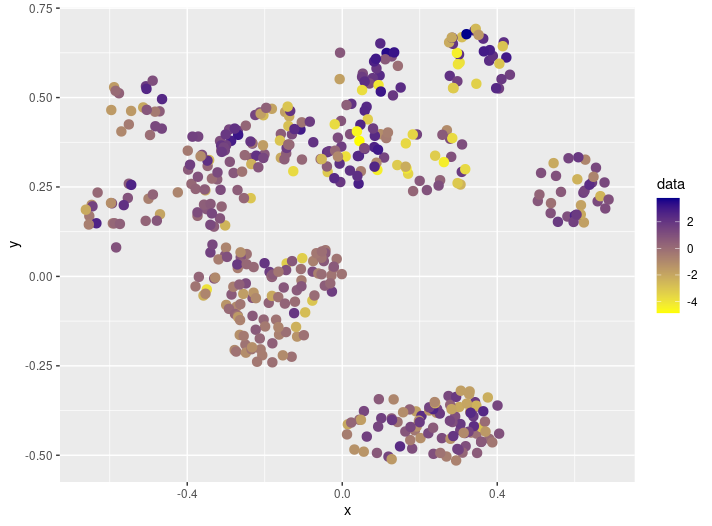}\par 
    \textbf{$t=10$}\par
    \includegraphics[width=0.6\linewidth]{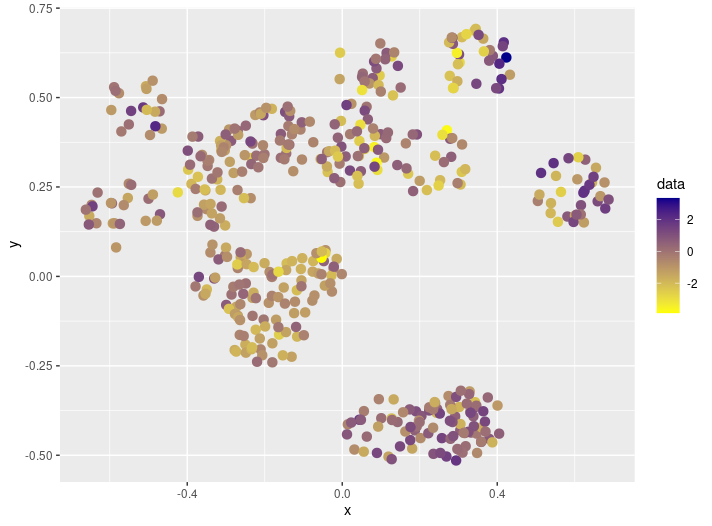}\par 
    \end{multicols}
\caption{Temporal evolution of the observations for structured selected locations. Process with only varying weight on the left and process with both varying weights and atoms on the right.}
\label{fig:simulations2_time}
\end{figure}

Figure \ref{fig:simulations2_time} shows the temporal evolution of the observations when locations are generated from a Thomas process. 

\newpage

\section*{Appendix D: An application to temperature data}

\begin{figure}
\begin{multicols}{4}
    \textbf{$1985$}\par
    \includegraphics[width=\linewidth]{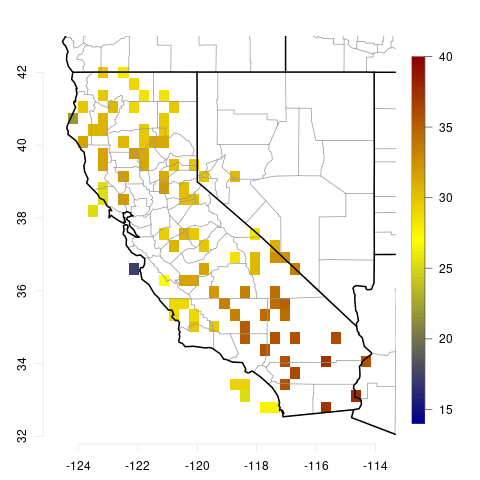}\par 
    \textbf{$1986$}\par
    \includegraphics[width=\linewidth]{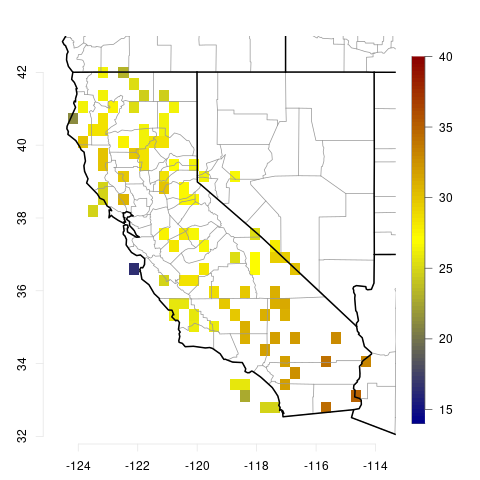}\par 
    \textbf{$1987$}\par
    \includegraphics[width=\linewidth]{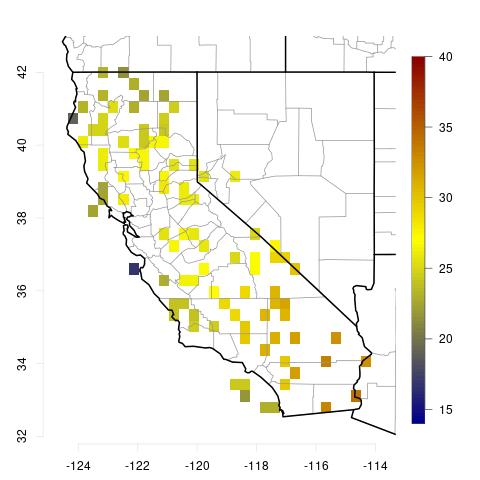}\par 
    \textbf{$1988$}\par
    \includegraphics[width=\linewidth]{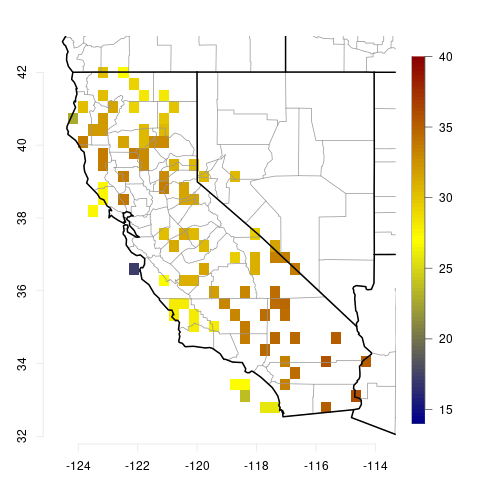}\par 
    \end{multicols}
\begin{multicols}{4}
    \textbf{$1989$}\par
    \includegraphics[width=\linewidth]{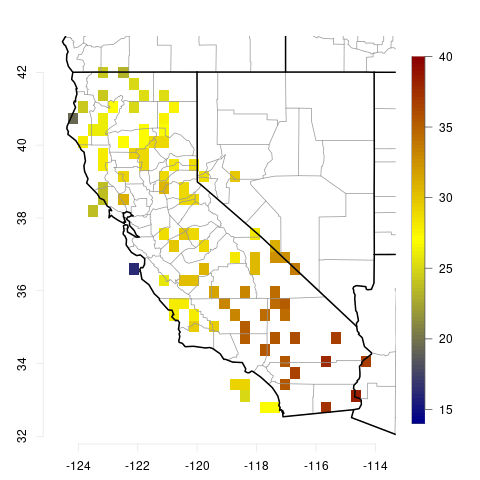}\par 
    \textbf{$1990$}\par
    \includegraphics[width=\linewidth]{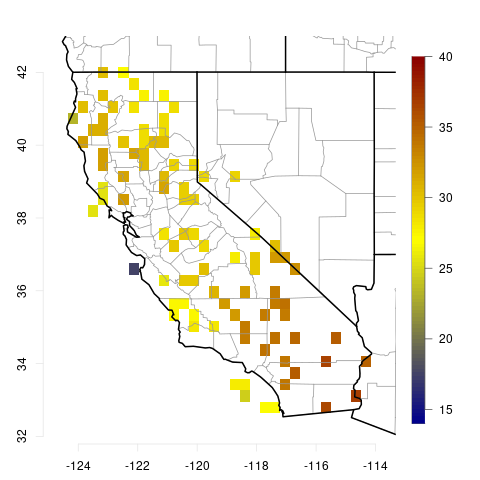}\par 
    \textbf{$1991$}\par
    \includegraphics[width=\linewidth]{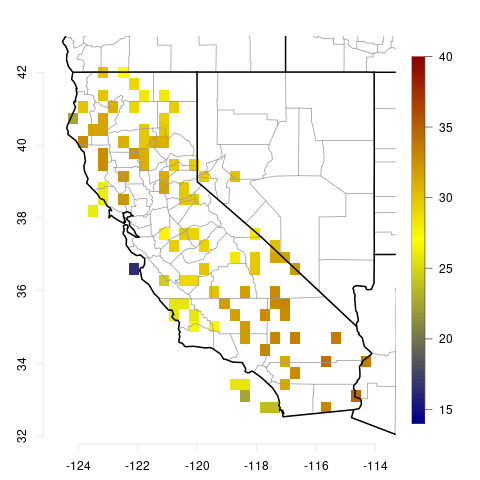}\par 
    \textbf{$1992$}\par
    \includegraphics[width=\linewidth]{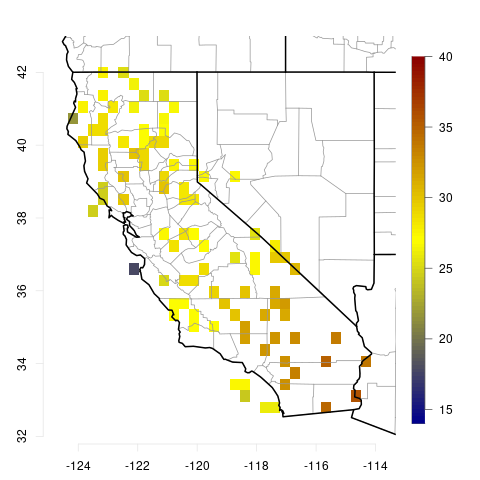}\par 
    \end{multicols}
\begin{multicols}{4}
    \textbf{$1993$}\par
    \includegraphics[width=\linewidth]{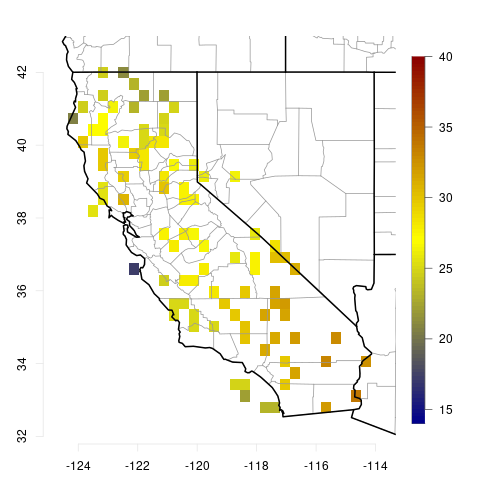}\par 
    \textbf{$1994$}\par
    \includegraphics[width=\linewidth]{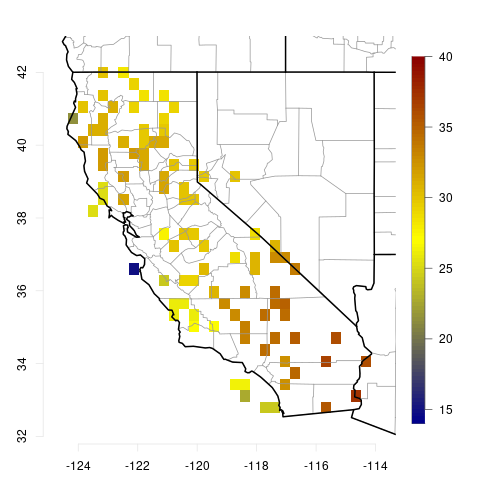}\par 
    \textbf{$1995$}\par
    \includegraphics[width=\linewidth]{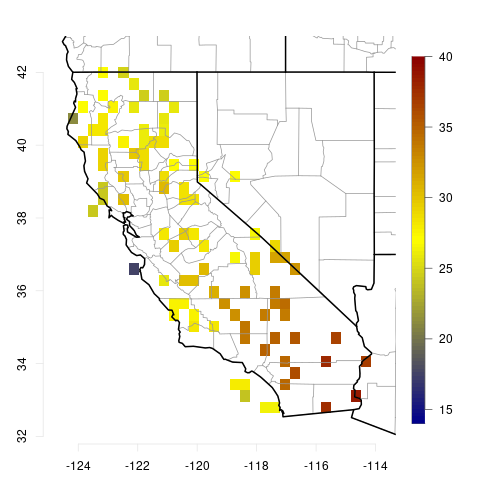}\par 
    \textbf{$1996$}\par
    \includegraphics[width=\linewidth]{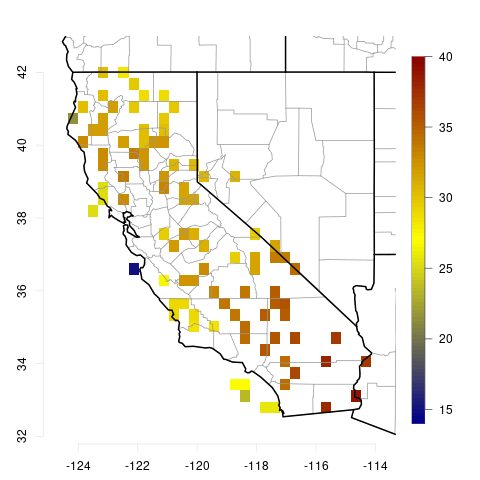}\par 
    \end{multicols}
\begin{multicols}{4}
    \textbf{$1997$}\par
    \includegraphics[width=\linewidth]{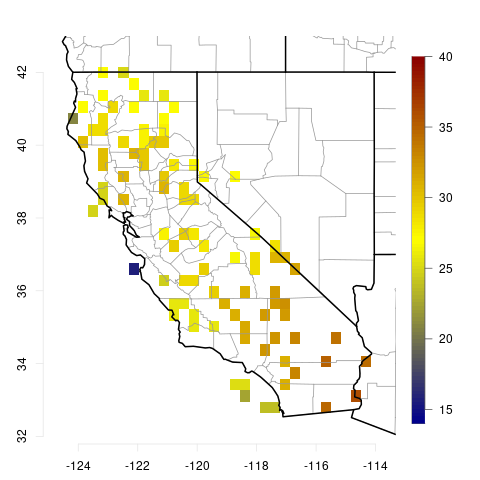}\par 
    \textbf{$1998$}\par
    \includegraphics[width=\linewidth]{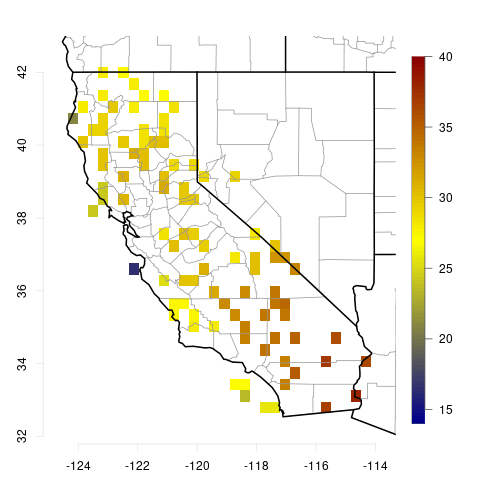}\par 
    \textbf{$1999$}\par
    \includegraphics[width=\linewidth]{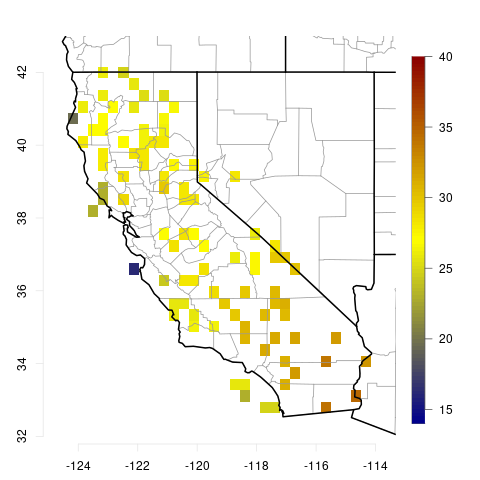}\par 
    \textbf{$2000$}\par
    \includegraphics[width=\linewidth]{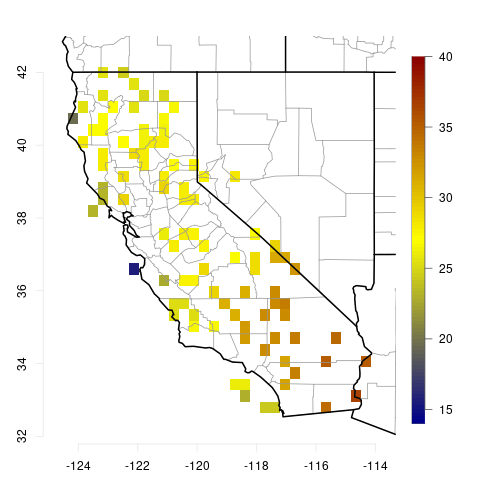}\par 
    \end{multicols}
\begin{multicols}{4}
    \textbf{$2001$}\par
    \includegraphics[width=\linewidth]{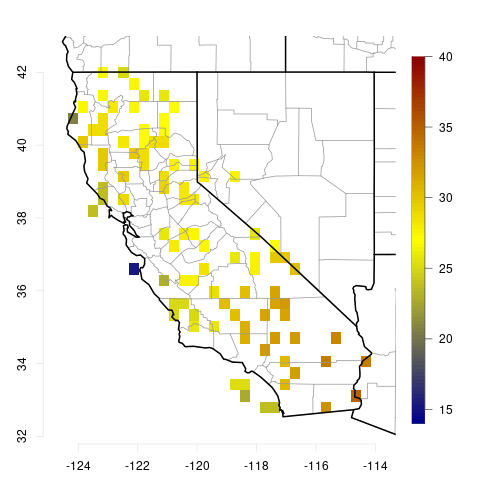}\par 
    \textbf{$2002$}\par
    \includegraphics[width=\linewidth]{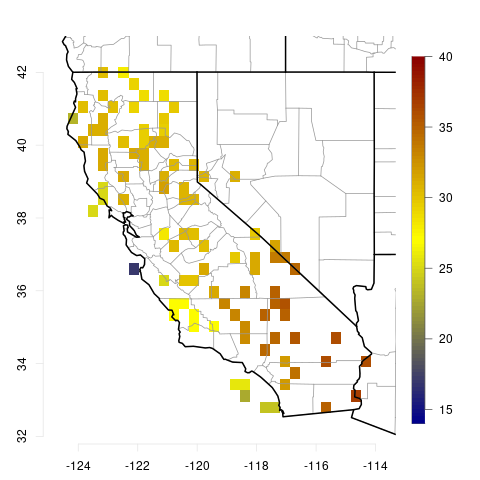}\par 
    \textbf{$2003$}\par
    \includegraphics[width=\linewidth]{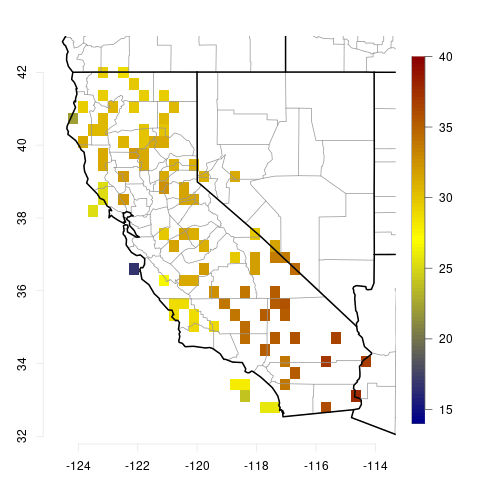}\par 
    \textbf{$2004$}\par
    \includegraphics[width=\linewidth]{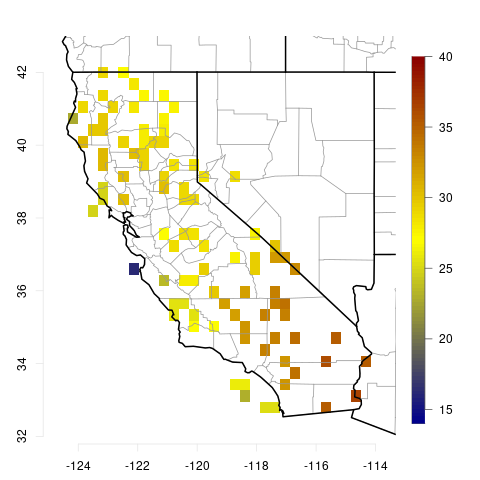}\par 
    \end{multicols}
\caption{Average maximum daily temperatures at the 100 locations in California from 1985 to 2004. Hotter locations are represented in red and cooler locations are represented in blue. Only July for each year is shown.}
\label{fig:temp_loc}
\end{figure}

The first dataset consists of daily averaged maximum temperatures at 10 meters for every month from 1985 to 2004, at 100 locations in California, USA. 
Data are gathered from the Atmospheric Science Data Center at NASA. For illustration, average temperatures for the month of July for each year are shown in Figure \ref{fig:temp_loc}, which also shows the locations of the monitoring stations. For the analysis, daily temperatures are used.

Four methods have been implemented, using daily temperatures from 1985 to 2003 as training dataset, and predicting the daily temperatures for year 2004: a spatio-temporal kriging model with coregionalisation (sp), a spatial stick-breaking (sSB) model with temporal varying atoms, and the two spatio-temporal stick-breaking models (stSB and stSB (VA)) proposed in this paper.

\begin{figure}
\centering
\begin{subfigure}{.33\textwidth}
  \centering
  \includegraphics[width=4.5cm,height=5cm]{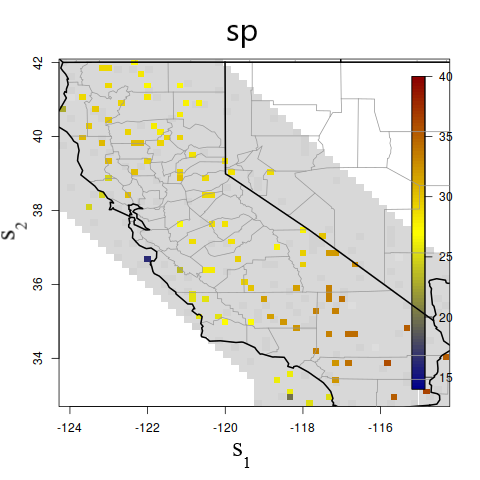}
\end{subfigure}%
\begin{subfigure}{.33\textwidth}
  \centering
  \includegraphics[width=4.5cm,height=5cm]{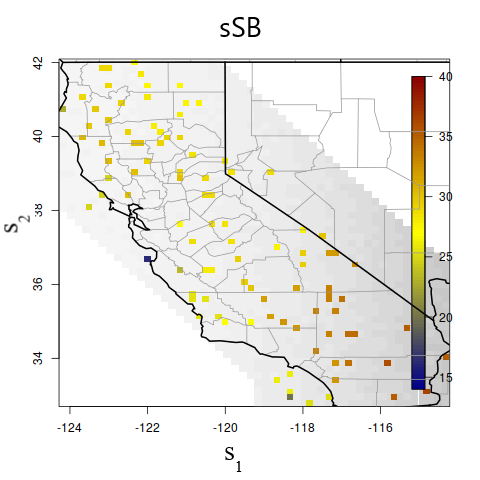}
\end{subfigure}
\begin{subfigure}{.33\textwidth}
  \centering
  \includegraphics[width=4.5cm,height=4.5cm]{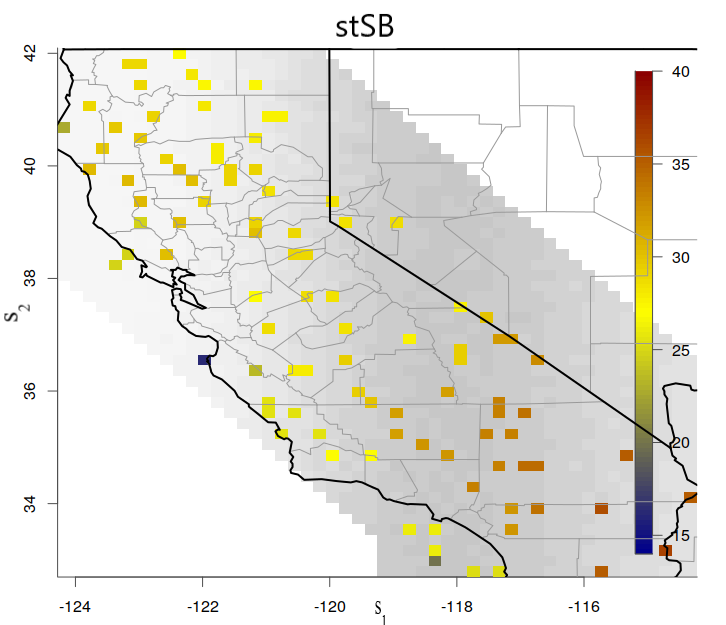}
\end{subfigure}
\caption{Average predictions for July 2004 given by sp (left), sSB (centre), and stSB (right). Darker gray is relative to hotter temperatures and lighter gray to cooler temperatures). The colored points represent the observations.}
\label{fig:rd_pred}
\end{figure}

Figure \ref{fig:rd_pred} shows the predicted means for July 2004 obtained with the three models, compared with the observed average temperatures. Similar comments can be applied to other months. The spatial model (left) oversmooths among the area of observation, with temperatures close to each other everywhere. The spatial stick-breaking (centre) shows a stronger spatial structure, however it tends to smooth locally: sudden and local changes over space and time are not well predicted; for example, the north-west area, where there is a small relatively hot area, is not well predicted. Figure \ref{fig:temp_loc} seems to suggest that this cluster is relatively recent (only years 2002-2003 show hotter temperatures there). On the contrary, the spatio-temporal stick-breaking model (right) is able to better incorporate heterogeneity, in particular to better adapt to recent changes: there are several local areas where the predicted temperature is better predicted than with other models. 

The better performance of stSB is reflected by the computation of the ESPE, across all days and all areas for year 2004: this is 103.66 for sp, 123.67 for sSB, and 40.53 for stSB. The stSB (VA) has also been implemented, but on a smaller training dataset, representing 50\% of the observations randomly selected; this has been done because of computational and memory reasons. The corresponding ESPE is 96.99; the model still performs better the sp and sSB in terms of predictive accuracy, however, it is trained on a smaller number of observations with an additional complexity in the number of parameters. 

Finally, it is possible to study the separability of the spatial and temporal component in the definition of the stick-breaking probabilities, by analysing the posterior distribution of the parameter $\lambda$, which is given a spike-and-slab prior, with one component concentrated around zero and one component given by a beta density $Be(a,b)$ with shape parameters $a=1$ and $b=1$. We have obtained that 100\% of the times $\lambda \neq 0$ a posteriori, which is against the assumption of separability, for both the single-atom model and the model with varying weights and atoms.  

\end{document}